\documentclass[11pt]{amsart}

\usepackage{amscd,amsmath,amsfonts,amssymb}

\newcommand{\la}{\lambda}
\newcommand{\al}{\alpha}
\newcommand{\be}{\beta}

\newcommand{\f}{\varphi}
\newcommand{\e}{\varepsilon}

\newcommand{\norm}[1]{\Vert #1\Vert}
\newcommand{\abs}[1]{\vert #1\vert}
\newcommand{\ov}{\overline}
\newcommand{\Ll}{\mathcal{L}}

\newcommand{\Ker}{\ensuremath\mathop\mathrm{Ker}}

\newcommand{\N}{\ensuremath\mathop\mathrm{Null}}

\newcommand{\hhat}[1]{\widehat #1}
\newcommand{\ra}{\rightarrow}
\newcommand{\dt}{d^\theta}

\newcommand{\CC}{\mathbb{C}}

\newcommand{\RR}{\mathbb{R}}
\newcommand{\ZZ}{\mathbb{Z}}
\newcommand{\NN}{\mathbb{N}}
\newcommand{\QQ}{\mathbb{Q}}

\newcommand{\Kmin}{\ensuremath{K_{\text{min}}}}
\newcommand{\Kminprime}{\ensuremath{K_{\text{min}}'}}

\newcommand{\Gammamin}{\ensuremath{\Gamma_{\text{min}}}}
\newcommand{\Gammaminprime}{\ensuremath{\Gamma_{\text{min}}'}}
\newcommand{\presmin}{\ensuremath{(\Kmin,\Gammamin)}}
\newcommand{\presminprime}{\ensuremath{(\Kminprime,\Gammaminprime)}}
\newcommand{\Kmax}{\ensuremath{\tilde{K}}}
\newcommand{\Kmaxprime}{\ensuremath{\tilde{K}'}}

\newcommand{\Gammamax}{\ensuremath{\tilde{\Gamma}}}
\newcommand{\Gammamaxprime}{\ensuremath{\tilde{\Gamma}'}}
\newcommand{\presmax}{\ensuremath{(\Kmax,\Gammamax)}}
\newcommand{\presmaxprime}{\ensuremath{(\Kmaxprime,\Gammamaxprime)}}

\newcommand{\Isom}[1]{\ensuremath{\text{\upshape \rmfamily Isom}(#1)}}
\newcommand{\map}[3]{\mbox{${#1}\colon{#2}\to{#3}$}}
\newcommand{\rid}{\ensuremath{/\!\!/}}

\numberwithin{equation}{section}

\newtheorem{te}{Theorem}[section]
\newtheorem{pr}{Proposition}[section]
\newtheorem{co}{Corollary}[section]

\theoremstyle{definition}
\newtheorem{de}{Definition}[section]
\newtheorem{re}{Remark}[section]
\newtheorem{ex}{Example}[section]
\newtheorem{cex}{Counterexample}[section]

\begin{document}

\title[Locally conformally K{\"a}hler manifolds]{Locally conformally
 K{\"a}hler manifolds.\\ A selection of results.}
\thanks{
 The author
thanks the Universit\`a della Basilicata (Potenza)  and Erwin
Schr\"odinger Institute (Vienna) for hospitality and financial
support during July, respectively August 2003, when part of this
paper was written.}

\author[Liviu Ornea]{Liviu Ornea}

 \address{\newline
\noindent University of Bucharest,\newline Faculty of
Mathematics,\newline 14 Academiei str.,\newline  70109 Bucharest,
Romania}

\email{ Liviu.Ornea@imar.ro, lornea@gta.math.unibuc.ro}

\begin{abstract}
I present a selection of results on locally conformally K\"ahler
geometry published after 1997. The proofs are mainly sketched,
some of them are even omitted. Several open problems are indicated in the end.
\end{abstract}
\maketitle



\tableofcontents

\section{Foreword}

Although known since the 50's, from P. Libermann's work, locally
conformally K{\"a}hler (l.c.K.) structures have been
systematically studied only starting with 1976, after the impulse given by  I.
Vaisman's work. Twen\-ty-two years after, the monograph \cite{drag}
gathered almost all known results in the subject. But already
during the printing of the book new results where announced  that
could not be reported. Somehow, I had the feeling that the book
came a little bit too early. On the other hand, I cherish the hope
that it gave an impetus for new research.

Soon after the book was published, P. Gauduchon and the present
author proved the existence of l.c.K. structures on all Hopf
surfaces, a construction later generalized in the paper 
\cite{ko} to produce Hopf manifolds in any dimension. Almost in
the same time, in his thesis, F.A. Belgun classified the compact
complex surfaces with l.c.K. metrics with parallel Lee form
(Vaisman metrics). In the same time, Y. Kamishima studied the
l.c.K. structures from the uniformization viewpoint in a series of
papers. He also initiated the study of the automorphism group of
l.c.K. manifolds. The properties of the Ricci-like curvature were
investigated by authors like B. Alexandrov, S. Ivanov; using
Bochner type formulas, they found some new topological properties
of l.c.K. manifolds with positive Weyl-Ricci tensor. With a 
completely different approach, M. Verbitsky recently produced similar
restrictions and a structure theorem for compact Einstein-Weyl
Vaisman manifolds, relating this geometry to the Sasakian one,
very much in the spirit of the characterization given in \cite{ko}
and \cite{gop}. The latter paper also discusses Hamiltonian fields
and extends the reduction procedure from symplectic and K{\"a}hler
geometry to l.c.K. structures. This study was further developed in
\cite{gopp}, where a new equivalent definition was proposed, in terms
  of presentations; it makes more clear the relation between K\"ahler and
 l.c.K. geometries and provides a new invariant for l.c.K. structures.
 The
structure of compact Vaisman manifolds is now  completely understood
in \cite{ov1} and, on the other hand, algebro-geometric techniques
were applied in order to derive vanishing theorems and a Kodaira
type immersion theorem for Vaisman manifolds, see \cite{ve},
\cite{ov2}. Moreover,  Riemannian properties of the l.c.K. 
metrics, in particular harmonicity of transversal holomorphic maps and 
harmonicity and minimality of vector fields and distributions 
 were also discussed in \cite{ba_dr} and \cite{orv}.
Very recently, also the Inoue and Inoue-Hirzebruch surfaces 
were generalized to higher dimension, in \cite{ot} 
(where, as a side result, the authors also disproved the 30 years 
lasting Vaisman conjecture stating that any compact 
l.c.K. manifold has an odd odd Betti number) and \cite{re}. 
Further on, considering in \cite{ov3} 
l.c.K. manifolds with potential, 
a class including the Vaisman one, 
the Hopf surfaces of K\"ahler rank $0$ were also 
generalized to higher dimension; any compact l.c.K. manifold with 
potential, of complex dimension at least $3$, can be embedded in the 
resulting linear Hopf manifold. Finally, using the Gauduchon metric on 
compact Vaisman manifolds, a stability theory was developed in 
\cite{ve04}, with particular striking results for stable bundles over 
diagonal Hopf manifolds. 

All this development is intimately connected with the late
achievements in Sasakian and 3-Sasakian geometry. One may argue
that the structure theorem in \cite{ov1} reduces Vaisman geometry, in the compact case, 
to Sasakian one. But on the other hand, also results in l.c.K. geometry served to obtain significant progress in Sasakian geometry: Belgun used his
classification of compact complex surfaces
admitting metrics with parallel Lee form to completely classify
Sasakian structures on compact 3-manifolds, cf.  \cite{be1}, \cite{be2},
 and the
embedding theorem for compact Vaisman manifolds produced, cf. \cite{ov2} and \cite{ov3},  a
CR-embedding of any compact Sasakian
manifold into a Sasakian weighted sphere. However, it is worth 
stressing that all these results hold only in the compact case.

Let me finally mention that much work was done in the late years
on understanding the locally conformally K\"ahler structures in
quaternionic geometry, but I shall not report here on this topic.
The interested reader may consult my former survey \cite{or} and
the more recent paper \cite{ve}.

Needless to say, the present survey does not aim to exhaustive
comprehension. For example, I have left the recent results on 
special submanifolds of l.c.K. manifolds (cf. \cite{barl}, \cite{bl_dr} \emph{et. al.}) 
for a future paper. This report  only reflects my knowledge, 
information, taste and power of
understanding.

\medskip

\noindent {\bf Acknowledgements.} I thank all the people who 
are interested in l.c.K. geometry, giving me the opportunity 
to write this text. 

I am particularly grateful to all my collaborators 
and friends who taught me a lot of mathematics. 

I  thank Alex Buium,  
Karl Oeljeklaus and Matei Toma who patiently explained me the paper \cite{ot} and 
Florin Belgun, Maurizio Parton, Misha Verbitsky and Victor Vuletescu for carefully 
reading a first draft of the paper, for correcting many of my mistakes and for suggesting important improvements.

\section{Review of old results and a new definition}

\subsection{Definitions} 
For this subsection, I mainly refer to \cite{drag} and the references
therein.

In the sequel  $(M,J,g)$ will be a connected Hermitian manifold of
complex dimension $n \geq 2$.  I shall denote by $\omega$ its
fundamental two-form given by $\omega(X,Y)=g(X,JY)$.

\begin{de}
$(M,J,g)$ is called {\em locally conformally   K{\"a}hler},  {\em
l.c.K.} for short,  if there
 exists an open cover $\mathcal{U}=\{U_{\alpha}\}$ such that each
 locally defined metric
$g_{\alpha}=e^{-f_\alpha}g_{|U_\alpha}$ is K{\"a}hler for some smooth
      function
       $f_\alpha$ on $U_\alpha$.
\end{de}

Equivalently,   $(M,J,g)$ is l.c.K. if and only if there exists a
{\em closed} one-form $\theta$ such that
\begin{equation}\label{1}
d\omega=\theta\wedge  \omega.
\end{equation}

\noindent Of course, locally, $\omega_{|U}=df_U$. Note also that, except on
complex surfaces, the equation \eqref{1} implies $d\theta =0$.

The globally defined one-form $\theta$ is called {\em the Lee form\footnote{After the name of Hwa-Cwung Lee, cf. \cite{lee}.}} and its metrically
equivalent (with respect to $g$) vector field $B=\theta^\sharp$ is called {\em the Lee
vector field}; the vector field $JB$ is called the {\em anti-Lee vector field}.

It can be easily seen that  $(M,J,g)$ is
l.c.K. if and
only if the following equation is satisfied for any $X,Y\in
\mathcal{X}(M)$:
\begin{equation}\label{nablaj}
(\nabla_XJ)Y=\frac 12
\{\theta(JY)X-\theta(Y)JX+g(X,Y)JB-\omega(X,Y)B\},
\end{equation}

\noindent where $\nabla$ is the Levi-Civita connection of $g$.
Note that the above equation shows that l.c.K. manifolds belong to
the
 class $W_4$ of the Gray-Hervella classification.

Another very useful characterization of l.c.K. manifolds is obtained
\emph{via} the universal Riemannian cover. Let $p:\tilde M\rightarrow
M$ be the covering map and denote also by $J$ the lifted complex structure. The lifted
metric $\tilde g$ is globally conformal with a K{\"a}hler metric $\tilde
h$, because $p^*\theta$
is exact. The fundamental group of $M$ will then act by holomorphic
conformal transformations with respect to $(J,\tilde h)$. But
conformal transformations of a symplectic form (in real dimension at
least $4$) are in fact homotheties. The converse is also true. But, in fact, one can see that it is not important to work with the universal cover, but with 
any cover on which the pull-back of the Lee form is exact. Summing up, the following equivalent definition holds: 

\begin{te}
A complex manifold $(M,J)$ is l.c.K. if and only if admits a covering 
endowed with a global  K{\"a}hler metric with respect to which the deck group 
of the covering 
acts by holomorphic homotheties.
\end{te}

I shall denote by $\mathcal{H}(\tilde M, \tilde
h,J)$ the above group of holomorphic homotheties and by
$\rho:\mathcal{H}(\tilde M, \tilde h,J)\rightarrow \RR_+$ the
homomorphism associating to each homothety its scale factor.

\subsection{Presentations}
For this subsection I refer to \cite{gopp} where all the definitions and results are to be found.

The above theorem suggests the following  definition (where by homothetic K\"ahler one understands theta the K\"ahler metric is fixed only up to homotheties):
\begin{de}
A pair $(K,\Gamma)$ is a {\em presentation} 
if $K$ is a homothetic K\"ahler manifold and $\Gamma$ a 
discrete Lie group of biholomorphic homotheties acting freely and properly
discontinously on $K$. If $M$ is l.c.K. and $M=K/\Gamma$ 
as l.c.K. manifolds, then $(K,\Gamma)$ is a presentation {\em of} $M$.
\end{de}

There can be more than one presentation for a given l.c.K. manifold. One cannot give a hierarchy, but can distinguish the extremes and, consequently, give a new definition of l.c.K. manifolds:
\begin{de} 
et $(K,\Gamma)$ be a presentation.
The associated {\em maximal} presentation is \presmax,
where $\Kmax$ is the homothetic K\"ahler universal covering of $K$ and 
$\Gammamax$ is the lifting of $\Gamma$ to $\Kmax$\footnote{See \cite{gopp} for a thorough discussion about the properties of lifted actions.}.
Then 
\[
\presmin=\left(\frac{\Kmax}{\Isom{\Kmax}\cap\Gammamax},\frac{\Gammamax}{\Isom{\Kmax}\cap\Gammamax}\right)
\]
is called the associated {\em minimal} 
presentation.

The presentations  
$(K,\Gamma)$ and  $(K',\Gamma')$ are equivalent 
if $\presmax=\presmaxprime$. The $=$ sign means that there exists an equivariant map,
namely, there exists
a pair of maps $(f,h)$,  where \map{f}{\Kmax}{\Kmaxprime} is a biholomorphic homothety, 
\map{h}{\Gammamax}{\Gammamaxprime} is an isomorphism and for any $\gamma\in\Gammamax$,
$x\in\Kmax$ we have
$f(\gamma x)=h(\gamma)f(x)$.  
Equivalently, we say that 
$(K,\Gamma)$ is equivalent to $(K',\Gamma')$
if  
$\presmin=\presminprime$. 

A l.c.K. manifold is an equivalence class $[(K,\Gamma)]$ of presentations.
\end{de}

Of course, the second part of the definition is motivated by the fact that the maximal and minimal presentations are uniquely associated (as $\Gammamax$, $\Gammamin$-spaces respectively) to a given l.c.K. structure.

Note that any minimal presentation is necessarily $(K,\ZZ)$.

In this context, I shall denote with $\rho_K$ the  scale homomorphism (the character)  associated to the presentation $(K,\Gamma)$, $\rho_K:\Gamma\rightarrow \RR_+$. Its image is a finitely generated subgroup (see also \cite{ov3}) and its rank provides a nice invariant of the l.c.K. structure, measuring the "trully conformal" part of $\Gamma$:

\begin{pr}\label{rang}
For any presentation $(K,\Gamma)$, the rank of the free abelian group 
$\rho_\Gamma(\Gamma)$  
depends only on the equivalence class $[(K,\Gamma)]$. Hence, one can speak about the rank $R_M$ of l.c.K. manifold. 
\end{pr}
The rank is bounded:
\begin{pr}\label{rM} 
\[
0\leq r_M\leq b_1(M)
\]
and $r_M=0$ if and only if $M$ is
g.c.K. In particular, if $b_1(M)=1$ then $r_M=1$.
\end{pr}
For example, the manifolds $M_s=(H^s\times \CC)\Gamma$ constructed in \cite{ot} (cf. Section \ref{matei}) all have rank $s$. 

\subsection{Vaisman manifolds}
A strictly smaller class of l.c.K. manifolds is the one formed by
those with parallel (with respect to the Levi Civita connection)
Lee form. I call  them {\em Vaisman manifolds}, as I. Vaisman was
the first to study them systematically (under the name of
generalized Hopf manifolds,
 \cite{vai}, a name which later proved to be inappropriate).
On such a  manifold, the length of the Lee vector field is
constant and I shall always assume it is nonzero. Hence, in what
follows, I shall normalize and consider that on a
 Vaisman
manifold $\norm{B}=\norm{JB}=1$.

Note that on a compact l.c.K. manifold, the metric with parallel Lee form, if
it exists, is unique up to homothety in its conformal class and
coincide with the Gauduchon (standard) metric, \cite{mpps}.

\begin{pr}\label{kill}
Let $(M,J,g)$ be a Vaisman manifold. Then
the Lee and anti-Lee vector fields commute $([B,JB]=0)$, are Killing
$(\Ll_Bg=\Ll_{JB}g=0)$ and holomorphic $(\Ll_BJ=\Ll_{JB}J=0)$.
Consequently, the distribution generated by $B$ and $JB$ is a
holomorphic Riemannian foliation.
\end{pr}

I denote by $\mathcal{F}$ the foliation generated by $B$ and $JB$.
Note also that the leaves of the  foliation generated by the
nullity of the Lee form carry an induced $\alpha$-Sasakian
structure (see \cite{blair} as concerns metric contact manifolds)
with $JB$ as characteristic (Reeb) vector field. If these foliations are
quasi-regular, then the leaf spaces are, respectively, a
K{\"a}hler and a Sasakian orbifold.

\medskip

The first known examples of l.c.K. manifolds which are not K{\"a}hler
(because they have $b_1=1$)
were the (diagonal) Hopf manifolds: $H_\lambda=(\CC^n\setminus \{0\}/\ZZ, \frac
{g_{flat}}{\abs{z}^2}, J_{can})$, where $\ZZ$ is generated by the
transformation $z_j\mapsto \lambda z_j$, with $\lambda\in \CC$,
$\abs{\lambda}\neq 0,1 $. The Lee form here is (when read on
$\CC^n\setminus \{0\}$): $\theta=-d\log \abs{z}^2$ and is seen to be
parallel. These manifolds are diffeomorphic with the product
$S^1\times S^{2n-1}$. What matters here is the Sasakian structure of
$S^{2n-1}$. Indeed it was then proved that the total space of a flat principal circle bundle over a
compact Sasakian manifold carries a Vaisman metric whose Lee form is
identified with the connection form of the bundle.

A Vaisman structure obtained as above is regular. Non-regular
examples can be obtained as suspensions over $S^1$ with fibre a
Sasakian manifold, \cite{go} and, for a general structure theorem,
\cite{ve}, \cite{ov1} (see below).

\begin{re} Particular examples of Vaisman structures are the locally conformally hyperk\"ahler ones. If compact, they have $b_1=1$, hence also rank $r_M=1$. We shall see later that, in fact, all compact Vaisman manifolds have rank $1$. On the other hand, examples of compact Vaisman manifolds with arbitrarily large $b_1$ can be obtained as induced Hopf bundles over curves of large genus in $\CC P^2$.
\end{re}

L.c.K. metrics without parallel Lee form were constructed by F.
Tricerri on some of the Inoue surfaces and by I. Vaisman on the
diagonal Hopf manifold.

\section{Locally conformally K\"ahler metrics on non-K{\"a}hler compact complex surfaces}

\subsection{Vaisman structures}
Locally, non-global, conformal K{\"a}hler metrics may exist only on
the surfaces in classes VI and VII of the Kodaira classification, \cite{kod}.

A generic Hopf surface is defined by Kodaira as a complex surface
whose universal cover is $\CC^2\setminus \{0\}$. Hopf surfaces with  fundamental
group isomorphic with $\ZZ$ are called primary.
Kodaira also proved that the fundamental group of any primary Hopf
surface can be realized as the cyclic group generated by the
transformation:
$$(z_1,z_2)\mapsto (\alpha z_1+\lambda z_2^m, \beta z_2),$$
where $m\in \NN$ and $\alpha, \beta, \lambda$ are complex numbers subject to the
conditions:
$$(\alpha-\beta ^m)\lambda=0, \quad \abs{\alpha}\geq\abs{\beta}>1.$$

A Vaisman metric on Hopf surfaces $H_{\alpha,\beta}$ with
$\lambda=0$ wass constructed in \cite{go}. It can be described as follows:

{\bf Step 1.}~ One starts with the canonical Sasakian structure $(
g_0, R_0)$ on $S^3$, where $g_0$ is the metric induced by the flat one
on $\CC^2$ and $R_0(z)=(iz_1, iz_2)$. With $k_1=\log\abs{\alpha}$,
$k_2=\log\abs{\beta}$, this one is deformed using the
function
$$\Delta(z)=\frac{2k_1\abs{z_1}^2+2k_2\abs{z_2}^2}{k_1+k_2},$$
namely one puts:
\begin{equation}
R_\Delta=R_0+\frac{2(k_1-k_2)}{k+1+k_2}(iz_1,-iz_2)
\end{equation}
and define a new metric $g_\Delta$ by the conditions:
\begin{enumerate}
\item $R_\Delta$ has $g_\Delta$-norm $1$.
\item $R_\Delta$ is orthogonal to the canonical contact distribution $\mathcal{D}$
  of $S^3$.
\item On $\mathcal{D}$, $g_\Delta=\Delta ^{-1} g_0$.
\end{enumerate}
Note that this is not a D-homothetic transformation in the sense of 
\cite{ta}. Still, the new structure $(g_\Delta, R_\Delta)$ is Sasakian
on $S^3$ (cf. \emph{loc. cit.} for a general criterion for such deformations to be Sasakian).

{\bf Step 2.} One uses the flow of $g_\Delta$-isometries
$$ \sigma_{\al, \be}((z,t))=\left(e ^{-i\mathfrak{Arg}\al}\cdot
  z_1,e ^{-i\mathfrak{Arg}\be}\cdot z_2), t\right)$$
to obtain a suspension over $S^1$ with fibre $S^3$. The total space
  will then be a Vaisman manifold diffeomorphic with $S^3\times S^1$.

To describe the underlying complex structure $J_{\al,\be}$, let $T$ be
the vector field tangent to $S^1$, viewed as a vector field on
$S^1\times S^3$, and let $E=(\bar z_2, -\bar z_1)$ (so that $E, iE$
generate $\mathcal{D}$). Also let
$$F(z)=\log \al\abs{z_1}^2 +\log \be\abs{z_2}^2.$$
Now $J_{\al,\be}$ is described by the following
table:
\begin{equation}
\begin{split}
J_{\al,\be}T&=\frac{1}{\mathfrak{Re}F} (-\mathfrak{Im}F
  T+\abs{F}^2R_0+i\bar F (\log \al-\log \be)z_1z_2E),\\
J_{\al,\be}R_0&=\frac{1}{\mathfrak{Re}F}(-T+\mathfrak{Im}FR_0+(\log
  \al-\log \be)z_1z_2E),\\
J_{\al,\be}E&=iE.
\end{split}
\end{equation}

Translated on $H_{\al,\be}$, by means of the diffeomorphism
$$F_{\al,\be}(z_1,z_2, t\mod\ZZ)=[\al^tz_1, \be ^tz_2],$$
this Vaisman structure corresponds to the potential
$$\Phi_{\alpha, \beta}(z_1,z_2)=e ^{\frac{(\log\abs{\alpha}+\log
    \abs{\beta})\tau}{2\pi}},$$
with $\tau$ given as (unique) solution as the equation:
$$\frac{\abs{z_1}^2}{e ^{\tau\log\abs{\alpha}\pi}}+
\frac{\abs{z_1}^2}{e ^{\tau\log\abs{\beta}\pi}}=1.$$
Namely:
\begin{te}\cite{go}
The metric:
$$g_{\alpha, \beta}=\frac{dd^c \Phi_{\alpha, \beta}}{\Phi_{\alpha,
    \beta}}(JX,Y)$$
is globally conformally K{\"a}hler on $\CC^{2}\setminus\{0\}$, has
    parallel Lee form and descends to a Vaisman metric on $H_{\alpha,\beta}$.
\end{te}

In fact, this way one gets a 1-parameter family of Vaisman metrics on
$H_{\alpha,\beta}$: indeed, for any $l\in \RR_+$, the potential
$\Phi_{\alpha, \beta}^l$ produces again a Vaisman metric (this amounts
to taking $lk_1, lk_2$ instead of $k_1, k_2$ above).

\medskip

This family of structures has been generalized in \cite{mau} by
including it in a larger family of l.c.K. structures (preserving the complex structure and
modifying only the metrics). Let
$h:\RR\rightarrow \RR_+$ be any function with period $2\pi$. Parton's
metrics are given by the following Hermitian matrix (written in the
complex basis $\{R_0,-E\}$):
\begin{equation}\label{h}
g_{\al,\be}^h=\left(
\begin{array}{ll}
\frac{\pi
  h}{(\mathfrak{Re}F)^2}+\frac{\abs{z_1}^2\abs{z_2}^2\log^2(\abs{\al}/\abs{\be}))}{(\mathfrak{Re}F)^3}
  &\frac{iz_1z_2\log(\abs{\al}/\abs{\be})}{(\mathfrak{Re}F)^2}\\
\\
\frac{\overline{iz_1z_2}\log(\abs{\al}/\abs{\be})}{(\mathfrak{Re}F)^2}&
\frac{1}{(\mathfrak{Re}F)}
\end{array}
\right)
\end{equation}
For a constant $h$, the family of Vaisman metrics in \cite{go} is
recovered for $l=2\pi h/(\log\abs{\al}+\log\abs{\be})$. This is not by
chance; indeed, by a direct computation one proves:
\begin{pr}\cite{mau} The metric $g_{\al,\be}^h$ is Vaisman if and only if
  $h=const.$
\end{pr}

 The
full list of compact complex surfaces which admit l.c.K. metrics with
parallel Lee form was given by Belgun:
\begin{te}\cite{belgun}\label{claspar}
A non-K{\"a}hler compact complex surface admits a Vaisman metric if and
only if it is: a properly
elliptic surface, or a Kodaira surface (primary or secondary), or a
Hopf
surface $H_{\al,\be}$.
\end{te}
Belgun's proof relies on the following criterion
for the existence of Vaisman metrics which encodes also the
construction in \cite{go}:
\begin{te}\cite{belgun}\label{crit}
Let $(M,J)$ be a compact complex surface with universal covering space
$(\tilde M,J)$. A Vaisman metric on $(M,J)$ is equivalent with the
following data:

$i)$~ A real holomorphic vector field $V$ without zeros on $M$.
 Denote with $\tilde V$ its lift to $\tilde M$ and with $\Phi$ its
complex flow.

$ii)$~ A group homomorphism $\tau: \pi_1(M)\cdot \Phi \ra
\RR$ such that
\begin{enumerate}
\item $\tau|_{\pi_1(M)}\neq 0$;
\item $\tau(\f_t^{\tilde V})=\e\neq 0$ and $\tau(\f_t^{J\tilde V})=  0$, where $\f_t^X$ denotes the flow of the vector field $X$.
\end{enumerate}

$iii)$~ A $3$-dimensional foliation of $\tilde M$ given by the level
hypersurfaces of a $\tau$ - equivariant function $f:\tilde M\ra \RR$
whose leaves carry a contact structure induced by  the distribution $\Omega$ of complex lines contained in
$\Ker df$.

In this case, a Vaisman metric on $M$ is given by the K{\"a}hler form
$\omega=k_1H^{-1}d(dH\circ J)$, where $H=e ^{k_2f}$ and $k_1,k_2\in
\RR^*$.
\end{te}
Here is the outline of the proof.

Let first
$g$ be a Vaisman metric  with (parallel) Lee form $\theta$ on $M$. By
the  de Rham
decomposition theorem,
$(\tilde M, \tilde g)= (\tilde S,g)\times (\RR, \text{can})$, where
$\tilde S$ is the universal cover of a leaf $S$ of the foliation
$\N\theta$ (note that $S$ carries an induces Sasakian structure which
underlying contact structure the statement refers to). Now $\tilde V$
is identified with $\norm{V}\partial_t$ and $f$ can be defined as the
projection of $\tilde M$ on the $\RR$ factor. To define $\tau$, one
first checks that for any $a\in \pi_1(M)$, the scalar $f((a\circ
\f^{\tilde V}_z)(x,0))$  does not depend on $x\in \tilde S$. Then one
sets:
$$\tau(a\circ\f^{\tilde V}_z)=  f((a\circ\f^{\tilde V}_z)(x,0)).$$
To show that $\tau$ is a group homomorphism, it is enough to see that
its restriction to $\pi_1(M)$ is such. This follows from the following
computation:
\begin{equation*}
\begin{split}
\tau(a_1a_2)&=f(a_1(a_2(x,0)))=f(a_1(\f_{\tau(a_2)}^{\tilde
  V}(x,0)))=f(\f_{\tau(a_2)}^{\tilde V}(a_1(x,0)))\\
&=f(\f_{\tau(a_2)}^{\tilde V}(\f_{\tau(a_1)}^{\tilde
  V}(x,0)))=\tau(a_1)+\tau(a_2).
\end{split}
\end{equation*}
The properties of $\tau$ are now straightforward ($\e=\norm{\tilde V}$).

The K{\"a}hler form is now defined for $k_1=-\frac 14$ and $k_2=2\e$.

For the converse, let $k_1=-\frac 14$ and $k_2=\pm 2df(\tilde V)$ with
the sign given by the one of $d(dH\circ J)(JY,Y)$, for $Y\in
\N\theta$. Then $\omega'= k_1H^{-1}d(dH\circ J)$ is clearly
l.c.K. with Lee form $\theta'=\frac{k_2}{2}df$ and Lee field
$\frac{k_2}{2}\tilde V$. This a unit Killing field (because
$\Ll_{\tilde V}f=0$ and $\Ll_{\tilde V}J=0$). Hence the symmetric part
of $\nabla \theta'$ is $0$. As the antisymmetric part of
$\nabla\theta'$is identified  with $d\theta'=0$, the proof is complete.

\begin{re}\label{cone}
For the first part of the proof one may also note the general fact
that
the cone $S\times\RR$ over a Sasakian manifold $(S,g)$ has the Vaisman
metric  $2e ^{-t}g_c$, with Lee form $-dt$, where $g_c=e ^t(dt^2+g)$ is the Riemannian cone metric.
\end{re}

Now the proof of Theorem \ref{claspar} proceeds by a case by case
analysis, using the above criterion and the following
consequence:
\begin{pr}\label{lift}
Let $W$ be the lift of the parallel Lee field of a Vaisman metric to
the universal cover. Then the orbits of $JW$ and of the cyclic group
generated by $a\circ\f^W_{-r}$ are relatively compact in
$\CC^2\setminus \{0\}$, where $a\in\pi_1(M)$ and $r=\tau(a)$.
\end{pr}

\subsection{Locally conformally K\"ahler metrics with non-parallel Lee form}
Other surfaces in the classes VI and VII -- but not all of them
-- may admit l.c.K. metrics with non-parallel Lee form.
The first examples (on surfaces which were not known to admit also
Vaisman structures) were constructed in \cite{tric} on the Inoue
surfaces $S_M$, $S_{n;p,q,r}^-$ and  $S_{n;p,q,u}^+$ with $u\in
\RR$. These metrics have harmonic (but non-parallel) Lee form and are
locally homogeneous.

A l.c.K. metric with non-parallel Lee form was
constructed in \cite{go} on the Hopf surface of K\"ahler rank 0
(with $\lambda\neq 0$). The construction uses a small deformation
of the structure on $H_{\al,\be}$, but the argument is specific to
this case: in general, the l.c.K. class is not closed under small
deformations \cite{belgun} (even if one can always deform a
Vaisman structure on a compact manifold to a quasi-regular one,
see below the proof of the embedding theorem in \cite{ov2}); see also \S \ref{lck_pot} for a subclass of l.c.K. manifolds which is stable to small deformations. 
Moreover, Belgun proved:
\begin{pr}\cite{belgun}
The Hopf surfaces of K\"ahler rank 0 do not admit any Vaisman
structure.
\end{pr}
The proof is again an application of the criterion in Theorem \ref{crit}. The
general form of a holomorphic vector field on $\CC^2\setminus \{0\}$
which descends on this Hopf surface can be found as:
$$W(z_1,z_2)=(mbz_1+cz_2^m)\partial_{z_1}+bz_2\partial_{z_2},\quad
b,c\in\CC.$$
Now one may compute the flows of $W$ and $JW$ and, supposing the
existence of a parallel Lee vector field (necessarily holomorphic),
apply Proposition
\ref{lift}: if the orbits of $JW$ are relatively compact, then $b\in
\RR$ and $c=0$. Then the second condition in the Proposition
\ref{lift} assures $\lambda=0$.

\medskip

Besides, Belgun  proved that the Inoue surfaces cannot admit l.c.K. 
metrics with parallel Lee form. And, most important, showing that compact non-K\"ahler surfaces are not necessarily of l.c.K. type, he proved:

\begin{te}\cite{belgun}
Inoue surfaces $S_{n;p,q,u}^+$ with $u\in
\CC\setminus \RR$ cannot admit any l.c.K. metric. 
\end{te}

This is in fact the only case where the method in \cite{tric} of constructing l.c.K. metrics on Inoue surfaces didn't work. The considered surface is a quotient of ${Sol_1'}^4$  by an integer lattice. Belgun's proof goes like this: he first shows that ${Sol_1'}^4$ admits a bi-invariant volume form; this is done by writing down explicitely the generators of the Lie algebra of ${Sol_1'}^4$ and checking that the multivector they define is $\mathrm{ad}$-invariant. He then proves that for $\mathfrak{Im} u\neq 0$, there is no left invariant l.c.K. metric on ${Sol_1'}^4$ (with another proof, this was also observed in \cite{vai1}). This follows from a very nice computation on the above determined generators of ${\mathfrak{sol}_1'}^4$ read as vector fields. Finally, he shows that if $S_{n;p,q,u}^+$ with $u\in
\CC\setminus \RR$ admits an l.c.K. metric, then, by averaging with respect to  to the  bi-invariant volume form found at the first step, it also admits a  ${Sol_1'}^4$-invariant l.c.K. metric, which is the desired contradiction. 

As $S_{n;p,q,u}^+$ with $u\in \CC\setminus \RR$ can be obtained by a small deformation of an Inoue surface which admits l.c.K. metric, this proves

\begin{co}\cite{belgun}
Unless the K\"ahler class, the l.c.K. class is not stable under small deformations.
\end{co}

On the other hand, as already mentioned, l.c.K. metrics with non-parallel Lee
form on some of the Inoue surfaces were constructed in \cite{tric}. I shall briefly recall this construction for $S_M$. Let
$H=\{w=w_1+iw_2\in\CC \mid w_2>0\}$ be the open half-plane and let $M=(a_{ij})\in
\mathrm{SL}(3,\ZZ)$ be a uni-modular matrix with one real eigenvalue $\al>1$ with
eigenvector $(a_1,a_2,a_3)$, and a non-real complex eigenvalue
$\beta$, with eigenvector $(b_1,b_2,b_3)$. Consider the following transformations:
\begin{equation*}
\begin{split}
(w,z)&\mapsto (\al w, \be z),\\
(w,z)&\mapsto(w+a_j, z+b_j)
\end{split}
\end{equation*}
They generate a group $\Gamma_M$ which 
acts on $H\times \CC$, the quotient being  a compact complex surface,
the Inoue surface $S_M$. The metric $g=w_2^{-2}dw\otimes d\bar
w+w_2dz\otimes d\bar z$ on $H\times \CC$ is globally conformal
K{\"a}hler with Lee form $\omega=d\log w_2$. Being compatible with the action
of $\Gamma_M$, it induces a l.c.K. metric on $S_M$.

This metric, a warped-product of the flat metric of $\CC$ and the $-1$-constant curvature metric of the Poincar\'e half-plane, has also interesting Riemannian properties. In fact, it was shown that its Lee and anti-Lee vector fields provide examples of harmonic vector fields (for definitions, see \emph{e.g.} \cite{BV00}, \cite{ggv}): 
\begin{pr}\cite{orv}\label{orv_harm} The Lee and anti-Lee vector fields of the Tricerri metric on an Inoue surface $S_M$ satisfy the following properties:

$i)$ they are harmonic and minimal;

$ii)$ the distribution locally generated by them (which is not a foliation in this case)  is harmonic and
determines a minimal immersion of $(S_M,g)$ into $(G^{or}_2(S_M),g^S)$.
\end{pr}

\begin{re} 
A l.c.K. structure on a principal $S^1$-bundle over a compact 3-dimensional 
solvable Lie group was constructed in \cite{acfm}. We reported it as such in 
\cite{drag}. In fact, as is was shown later on in \cite{ka4}, this structure is homothetically holomorphic with the l.c.K. structure on the Inoue surface of type $Sol^4_1$ constructed by Tricerri.
\end{re}

Recently, in \cite{fp}, a new construction was announced of anti--self--dual Hermitian metrics on some compact surfaces in class VII with $b_2>0$.  Such metrics are automatically l.c.K. by a result in \cite{bo}.

\subsection{A  characterization in terms of Dolbeault operator}
L.c.K. manifolds do not usually appear among the limiting cases of
inequalities involving the eigenvalues of Dirac-type operators.
One motivation could be that l.c.K. geometry has not a holonomy
formulation. However, such a characterization is available for
some l.c.K. surfaces in terms of the eigenvalues of the Dolbeault
operator (which, on K\"ahler manifolds, acts essentially as the
Dirac one). The proof is too technical to be reported here, I only
state the result:
\begin{te}\cite{agi}
Let $M$ be a compact Hermitian spin surface of positive conformal
scalar curvature $k$. Then the first eigenvalue of the Dolbeault
operator satisfies the inequality:
$$\lambda^2\geq \frac 12\, \mathrm{inf} k.$$
In the limiting case $k=const.$ and $M$ is l.c.K.
\end{te}

\section{A generalization of the Inoue surface $S_M$}\label{matei}

I shall outline here the construction in \cite{ot}. As far as I know, it is for the first time that (algebraic) number theory is used to construct examples of l.c.K. structures. It is remarkable that their examples also lead to disproving an almost thirty years old conjecture by Vaisman. 

Let $K$ be an algebraic number field of degree $n:=(K:\QQ)$. Let then $\sigma_1,\ldots, \sigma_s$ (resp. $\sigma_{s+1},\ldots,\sigma_n$) be the real (resp. complex) embeddings of $K$ into $\CC$, with $\sigma_{s+i}=\bar\sigma_{s+i+t}$, for $1\leq i\leq t$. Let  $\mathcal{O}_K$ be the ring of algebraic integers of $K$. Note that for any $s,t \in \NN$, there exist algebraic number fields with precisely $s$ real and $2t$ complex embeddings.

Using the embeddings $\sigma_i$, $K$ can be embedded in $\CC^m$, $m=s+t$, by 
$$\sigma:K\rightarrow \CC^m, \quad \sigma(a)=(\sigma_1(a),\ldots,
\sigma_m(a)).$$
This embedding extends to $\mathcal{O}_K$ and $\sigma(\mathcal{O}_K)$ is a lattice of rank $n$ in $\CC^m$, see, for example, \cite[p. 95 ff.]{bs}. This gives rise to a 
properly discontinuous action of $\mathcal{O}_K$ on $\CC^m$. On the other hand, $K$ itself acts on $\CC^m$ by
$$(a,z)\mapsto (\sigma_1(a)z_1,\ldots, \sigma_m(a)z_m).$$
Note that if $a\in \mathcal{O}_K$, $a\sigma(\mathcal{O}_K)\subseteq \sigma(\mathcal{O}_K)$. Let now $\mathcal{O}_K^*$ be the group of units in $\mathcal{O}_K$ and set
$$\mathcal{O}_K^{*,+}=\{a\in \mathcal{O}_K^* \mid \sigma_i(a)>0, \, 1\leq i\leq s\}.$$
The only torsion elements in the ring $\mathcal{O}_K^*$ are $\pm 1$, hence the Dirichlet units theorem asserts the existence of a free abelian group $G$ of rank $m-1$ such that $\mathcal{O}_K^* =G\cup (-G)$. Choose $G$ in such a  away that it contains $\mathcal{O}_K^{*,+}$ (with finite index). Now  $\mathcal{O}_K^{*,+}$acts multiplicatively on $\CC^m$ and, taking into account also the above additive   action, one obtains a free action of the semi-direct product $\mathcal{O}_K^{*,+}\ltimes \mathcal{O}_K^*$ on $\CC^m$ which leaves invariant $H^s\times \CC^t$ (as above, $H$ is the open upper half-plane in $\CC$). Again by Dirichlet units theorem the authors show that it is possible to choose a subgroup $U$ of $\mathcal{O}_K^{*,+}$ such that the action of $U\ltimes \mathcal{O}_K$ on $H^s\times \CC^t$ be properly discontinuous and co-compact. Such a subgroup $U$ is called \emph{admissible} for $K$. The quotient 
$$X(K,U):=(H^s\times \CC^t)/(U\ltimes \mathcal{O}_K)$$
is then shown to be a $m$-dimensional compact complex (affine) manifold, differentiably a fiber bundle over $(S^1)^s$ with fiber $(S^1)^n$.

Observe that for $s=t=1$ and $U=\mathcal{O}_K^{*,+}$, $X(K,U)$ reduces to an Inoue surface $S_M$ as described in the preceding paragraph.

\begin{te}\cite{ot}

$i)$ For $t=1$, $X(K,U)$ admits locally conformally K\"ahler metrics.

$ii)$ For $t>1$ and $s=1$, there is no locally conformally K\"ahler metric on $X(K,U)$.
\end{te}
Indeed, 
$$\f:H^s\times \CC\rightarrow \RR, \quad \f=\frac{1}{\Pi_{j=1}^s(i(z_j-\bar z_j))}+\abs{z_m}^2$$
is a K\"ahler potential on whose associated 2-form $i\partial\bar\partial \f$ the deck group acts by linear holomorphic homotheties.

A particular class of manifolds $X(K,U)$ is that of \emph{simple type}, when $U$ is not contained in $\ZZ$ and its action on  $\mathcal{O}_K$ does not admit a proper non-trivial invariant submodule of lower rank (which, as the authors show, is equivalent to the assumption that there is no proper intermediate field extension $\QQ\subset K'\subset K$ with $U\subset \mathcal{O}_{K'}$). The information about the topology of  $X(K,U)$ is gathered (I report only what is relevant to l.c.K. geometry, but the whole computation is beautiful) in the following
\begin{te}\cite{ot}

$i)$  $b_1(X(K,U))=s$, hence for $s=2p$ $X(K,U)$ cannot be Vaisman.

$ii)$ If $X(K,U)$ is of simple type, then 
$b_2(X(K,U))=\genfrac{(}{)}{0pt}{0}{s}{2}$.

$iii)$ The tangent bundle $T{X(K,U)}$ is flat and $\dim H^1(X(K,U), \mathcal{O}_{X(K,U)})\geq s$. In particular,  $X(K,U)$ are non-K\"ahler.
\end{te}

Now, for $s=2$ and $t=1$, the six-dimensional $X(K,U)$ is of simple type, hence has the following Betti numbers: $b_0=b_6=1$, $b_1=b_5=2$, $b_2=b_4=1$, by $i)$ and $ii)$,  $b_3=0$ by $iii)$. This proves:

\begin{te}\cite{ot}
Vaisman's conjecture claiming that \emph{a compact locally conformally K\"ahler, non-K\"ahler  manifold must have an odd odd Betti number} is false.
\end{te} 

This ends the story of a conjecture that most of the interested people believed true, especially because it holds for Vaisman manifolds ($b_1$ is odd). 

\begin{re}
No general classification of l.c.K. structures is available for the moment for $3$-dimensional complex manifolds. Only recently, in \cite{ugarte}, l.c.K. structure were classified on \emph{nilmanifolds}. Using representation theory, the author first determines the nilpotent Lie algebras that can admit l.c.K. structures: these are $\mathfrak{h}_1$ and $\mathfrak{h}_3$. He then proves that \emph{any l.c.K. structure on the product of the circle with a compact quotient of the Heisenberg group is Vaisman} and, on the other hand, \emph{compact complex parallelizable nilmanifolds which are not tori do not admit l.c.K. metrics.}
\end{re}

\section{The automorphism group of a locally conformally K\"ahler manifold}
By definition, a diffeomorphism of a l.c.K. manifold is a 
(l.c.K.)-automorphism if it is biholomorphic and preserves the
conformal class:
$$\mathrm{Aut}(M)=\{f\in \mathrm{Diff}(M) \mid f_*J=Jf_*, \,
f^*g\in [g]\}.$$
\begin{pr}\cite{ka2} The automorphism group of a compact l.c.K.
manifold is a compact Lie group.
\end{pr}
In fact, it is a Lie group as it is closed in the group of all
conformal transformations of $(M, [g])$ and it is compact as a
consequence of the Obata-Ferrand theorem.

This follows also from a more general argument. As I already
noted, on a compact Vaisman manifold the metric with parallel Lee
form coincides (up to homothety) with the Gauduchon metric. But
the Gauduchon metric exists on any compact l.c.K. manifold. Hence one may
thus think about the relation between $\mathrm{Aut}(M)$ and the
isometry group of the latter. And indeed:
\begin{te}\cite{mpps} The automorphism group of a compact l.c.K.
manifold coincides with the isometry group of the Gauduchon
metric.
\end{te}
This implies that, when working on compact Vaisman manifolds,
considering isometric actions with respect to the Vaisman metric
instead of conformal actions does not represent a loss of
generality (see Section \ref{red}).
\medskip

Vaisman structures can be characterized among the l.c.K. ones (but
only on  compact spaces) by means of group actions:
\begin{te}\cite{ko}\label{flux}
Let $(M,g,J)$ be a compact l.c.K. manifold. Then
$[g]$ contains  a metric with parallel Lee form if and only if
$\mathrm{Aut}(M)$ contains a complex $1$-dimensional Lie subgroup.
\end{te}
Clearly, from Proposition \ref{kill}, the condition is necessary.
The proof of the sufficiency is based on exhibiting the universal
cover $\tilde M$ of $M$ as a Riemannian cone (endowed with the associated K\"ahler structure) over a Sasakian manifold, then
showing that $\pi_1(M)$ acts by holomorphic homotheties with
respect to the K\"ahler metric. The needed Sasakian manifold is
obtained as follows. Let $T$ be the complex  1-dimensional Lie
subgroup in the statement. One first shows that a lift $\tilde T$
of $T$ cannot contain only isometries, hence $\rho(T)=\RR_+$.
It is then possible to pick
generators $\xi$, $J\xi$ such that the flow $\f^\xi_t$ act by
conformal maps (and is isomorphic to $\RR$) and the flow
$\f^{J\xi}_t$ act by isometries. Next, using the compactness of
$M$, one sees that this action of $\RR$ on $\tilde  M$ is free and
proper, in particular, $\xi$ does not vanish. Finally, the 1-level
set $W$ of the squared norm  of $\xi$ (with respect to the K\"ahler
metric of $\tilde M$) is shown to have an induced Sasakian
structure.

\medskip

On the other hand, it is known that the orbit map $ev$ of
any action of a torus $T^2$, $ev_x(t)=t\cdot x$,
induces a map $ev_*:\ZZ^2\longrightarrow H_1(M,\ZZ)$ in homology. If $M$ is
compact K\"ahler and the action is holomorphic, then  $ev_*$ is
injective. By contrast, Kamishima proves:
\begin{te}\cite{ka3}
Let $M$ be a compact, non-K\"ahler l.c.K. manifold
of complex dimension at least $2$. If $\mathrm{Aut}(M)$ contains a
complex torus $T^1_\CC$, then the induced action in homology has rank
$1$.
\end{te}
Note that some elliptic surfaces, in particular the diagonal Hopf
surfaces, more generally, the regular compact Vaisman manifolds,
 do admit such actions, generated by the Lee and anti-Lee vector
 fields.

\section{ Geometry and topology of Vaisman manifolds}
\subsection{A general example}\label{genex}
Examples of (compact), non-K{\"a}hler, l.c.K. manifolds are now
abundant. More precisely, there are many examples of compact
Vaisman manifolds. It was known that the total space of a flat
principal circle bundle over a compact Sasakian manifold carries a
Vaisman metric whose Lee form is identified with the connection
form of the bundle. I here present such a concrete example (see \cite{ko} for the general case
and \cite{go} for the surface case) which
plays in Vaisman geometry the role of the projective space in
K\"ahler geometry.

Let $a_j$ be real numbers such that $0<a_1\leq\ldots\leq a_n$. Deform the standard contact form of the sphere $S^{2n-1}$, $\eta=\sum (x_idy_i-y_ydx_i)$, to $\eta_A=\frac{1}{\sum a_i\mid z_i\mid^2}\eta$ which underlies the same contact structure, but whose Reeb field is $R_A=\sum a_i(x_i\partial y_i-y_i\partial x_i)$. Accordingly, define a metric $g_A$ on $S^{2n-1}$ by letting it be equal to the round one on $\Ker \eta_A$ and by requesting that $R_A$ be unitary and orthogonal to $\Ker \eta_A$. Now $(S^{2n-1}, \eta_A, g_A)$ is Sasakian and will be denoted $S_A^{2n-1}$. It is called {\em the weighted sphere}\footnote{This deformed Sasakian structure of the sphere has been recently encountered in a different context, related to the weighted projective spaces, cf. \cite{dg}}. The cone $\RR\times S^{2n-1}_A$ is then K\"ahler, \emph{e.g.} 
according to \cite{bgm}, the K\"ahler form being identified with $\omega_A=d(e^t\eta_A)$. The underlying complex manifold can be shown to be biholomorphic with $\CC^n\setminus\{0\}$ (with the standard complex structure) by means of the map $(t, (z_j))\mapsto (e^{-a_jt}z_j)$. Now pick a group $\Gamma\subset \RR\times \mathrm{PSH}(S^{2n-1}_A)$ which acts properly discontinuously by holomorphic homotheties with respect to $\omega_A$. The quotient space will be locally conformally K\"ahler with Lee form identified with $-dt$, hence parallel. Equivalently, we may take the corresponding action of $\Gamma$ on $\CC^n\setminus\{0\}$ and ``read'' this Vaisman structure on a quotient of  $\CC^n\setminus\{0\}$.  In particular, let 
$(c_1,\ldots,c_n)\in ({S^1})^n$, $s\neq 0$  and consider  an action of $\ZZ$ on $\CC^n\setminus\{0\}$ generated by: $(z_1,\ldots,z_n)\mapsto (e^{a_1s}c_1z_1,\ldots,e^{a_ns}c_nz_n)$. This group satisfies the above construction, hence the quotient $S^1\times S^{2n-1}_A$, isometric with $H_\Lambda:\CC^n\setminus\{0\}/\ZZ$, with $\la_i=e^{a_is}c_i$, is a Vaisman manifold that was called a (general, or non-standard, or diagonal) Hopf manifold. Reversing the construction, we see that for every complex numbers $\la_j$ satisfying $1<\mid \la_1\mid\leq \cdots\leq\mid\lambda_n\mid$, there exists a corresponding Hopf manifold.  In the simplest case, when $\lambda_i=1/2$, one
recovers the standard Hopf manifold with l.c.K. metric (read on
$\CC^n\setminus\{0\}$)
 $g_0=(\sum \abs{z_i}^2)^{-1}\sum dz_i\otimes d\bar z_i$ and Lee
form $\omega_0=-(\sum \abs{z_i}^2)^{-1}\sum (\bar z_i dz_i + z_i
d\bar z_i)$; here the Lee field is the one tangent to the $S^1$
factor.

The construction of the above example obeys a general pattern: one
starts with a Sasakian manifold, constructs the Riemannian cone over it, this has a K\"ahler structure, 
then selects a group which acts by holomorphic homotheties with
respect to this K\"ahler structure. The quotient is then a Vaisman
manifold. As will be clear later,  this pattern in fact exhausts
all compact Vaisman manifolds.

Still another way of looking at the above diagonal Hopf manifold, which will play an important role in \S \ref{lck_pot}, is from the 
viewpoint of the potential existing on the universal cover 
$\CC^n\setminus\{0\}$, cf. \cite{ve04}. In fact, if, for simplicity, we let $\Gamma$ be generated by the diagonal operator $A$ with eigenvalues $\al_i$, then the K\"ahler potential $$\f(z_1,\ldots,z_n)=\sum \abs{z_i}^{\be_i}, \quad \be_i=\log_{\abs{\al_i}^{-1}}C, \quad C=const.>1,$$
is acted on by  $A$ as follows: $A^*\f=C^{-1}\f$. Hence the associated K\"ahler form $\omega_K$ will satisfy the same relation: $A^*\omega_K=C^{-1}\omega_K$. This proves that the quotient $H_A=(\CC^n\setminus\{0\})/\langle A\rangle$
is l.c.K. Now Theorem \ref{flux} assures that it is Vaisman with respect to the metric $\frac{\omega_K}{\f}$ (remember also the construction in \cite{go}). It can be shown (cf. \cite{ve04}) that this is the Gauduchon metric of its conformal class and that its Lee field is given by the formula
$$\theta^\sharp=-\sum z_i\log \abs{\al_i}\partial z_i.$$

The diagonal Hopf manifold may be identified among the compact Vaisman manifolds by
means of the existence of a non-compact flow $\CC^*$ consisting of a kind of
transformations that I now describe. Let $\{\theta^1,\ldots,\theta^{n-1},\bar
\theta^1,\ldots,\bar\theta^{n-1}\}$ be complex 1-forms that together with $\theta$ and
$J\theta$ determine a coframe field adapted to the  l.c.K. structure of $M$. A diffeomorphism
$f$ of $M$ is called a \emph{Lie-Cauchy-Riemann transformation} (cf. \cite{ko}) if it
satisfies:
\begin{equation*}
\begin{split}
f^*\theta&=\theta,\\
 f^*(J\theta)&=\lambda J\theta,\\
f^*\theta^{\al}&=\sqrt \lambda\cdot\theta^\be U^{\al}_{\be}+
J\theta\cdot v^{\al},\\
f^*{\bar \theta}^{\al}&=\sqrt \lambda\cdot
{\bar\theta}^\be\ov{U}^\al_\be+J\theta\cdot\ov{v}^{\al},
\end{split}
\end{equation*}
for some smooth functions $\lambda\in \RR_+$,
$v^\al\in \CC$, $U^\al_\be\in{\rm U}(n-1)$. The totality of these
transformations form a group which contain the holomorphic isometries and
preserve the specific $G$-structure of l.c.K. manifolds. The announced result
is:
\begin{te}\cite{ko}
Let $(M,g,J)$ be a compact, connected, l.c.K. non-K\"ahler  manifold 
with  parallel Lee form $\theta$.
Suppose that $M$ admits a closed subgroup $\CC^*=S^1\times \RR_+$
of Lee-Cauchy-Riemann transformations whose $S^1$
subgroup induces the Lee field $\theta^\sharp$.
Then $M$ is holomorphically isometric, up to scalar multiple of the metric,
to a Hopf  manifold $H_\Lambda$.
\end{te}

\subsection{The weight bundle}\label{weight}
 Let now $L\longrightarrow M$ be the weight bundle canonically
associated to a l.c.K. manifold (see \cite{g_crelle}, \cite{cp}).
It is a real line bundle associated to the representation
$\mathrm{GL}(2n,\RR)\ni A\mapsto \mid \det A\mid^{\frac{1}{2n}}$. The Lee
form can be interpreted as a connection one-form in $L$,
associated to the Weyl connection determined on $M$ by the the
l.c.K. metric and the Lee form. $L$ is thus flat and one may speak
about its monodromy. Nothing can be said in general  about the
monodromy of $L$. However, for Vaisman manifolds, the following
result is true and it gives the key for the structure of such compact
manifolds:
\begin{te}\cite{ov1}\label{mono}
Let $M$ be a compact Vaisman manifold. Then the monodromy of $L$
is isomorphic to $\ZZ$.
\end{te}

For the proof, let $\Gamma$ be the monodromy (it is a subgroup of
$\RR_+$) and let $\hhat M\longrightarrow M$ be the associated
covering. Let $\mathrm{Aut}(\hhat M,M)$ be the group of conformal
automorphisms of $\hhat M$ which are lifts of a l.c.K. automorphism
of $M$. One may compose the scale homomorphism
$\rho:\mathcal{H}(\hhat M)\longrightarrow \RR_+$ with the natural
forgetful homomorphism $\Phi:\mathrm{Aut}(\hhat M,M)\longrightarrow
\mathrm{Aut}(M)$ and obtain a homomorphism
$$\mathrm{Aut}(\hhat M,M)\stackrel{\Phi\times\rho}{\longrightarrow}
\mathrm{Aut}(M)\times \RR_+$$ which proves to be injective. One
looks now to the Lie group $G$ generated by the Lee flow, which is
isomorphic to a torus ${(S^1)}^k$ (see above), and to its
counter-image in $\mathrm{Aut}(\hhat M,M)$, $\hhat G=\Phi^{-1}(G)$.
Observing that $\Gamma=\ker\Phi\subset\hhat G$, one is left to
prove that $\hhat G$ is isomorphic to $\RR\times {(S^1)}^{k-1}$.
This is done in two steps. One first shows that the connected
component $\hhat G_0$ of $\hhat G$ is non-compact, hence isomorphic
to $\RR\times {(S^1)}^{k-1}$. The second step here is to show that
$\hhat G$ is connected by showing that the quotient $H=\hhat G/G_0$
is trivial. To this end, one first shows that $H$ is finite, by
proving that each connected component of $\hhat G$ meets the
compact $\rho^{-1}$. Second, one notes that from the structure
theorem for abelian groups, $\hhat G\cong \RR\times
{(S^1)}^{k-1}\times H$, hence $H$ is embedded in $\hhat G$. As
such, it can be seen that $H$ is in fact embedded in $\Gamma$
which, as monodromy of a real line bundle, is included in $\RR$,
thus cannot have non-trivial compact subgroups.
\medskip

 The breakthrough towards obtaining new significant
results in l.c.K. geometry was the use of algebro-geometric
techniques by M. Verbitsky. The key tool of his approach is the
use of the 2-form
$$\omega_0:=d^c\theta=\frac{\partial -
\overline{\partial}}{\sqrt{-1}}\theta=-dJ\theta.$$ Its main
property is the following:
\begin{pr}\cite{ve}
Let $M$ be a Vaisman manifold. Then all eigenvalues of $\omega_0$
are positive, except the one corresponding to the Lee field which
is zero.
\end{pr}
On the other hand, it can easily be seen that on the orthogonal
transverse of the canonical foliation $\mathcal{F}$, $\omega_0$ is
closed, hence plays the role of a K\"ahler form. This implies
that, if the Vaisman manifold is quasi-regular, $\omega_0$
projects on the K\"ahler form of the base orbifold:
\begin{pr}\cite{ve}
Let $M$ be a compact, quasi-regular  Vaisman manifold. Let $p$
denote the canonical projection over the K\"ahler orbifold $(Q,
\omega_Q)$. Then $p^*\omega_Q=\omega_0$.
\end{pr}

The form $\omega_0$ can be interpreted as the curvature of the
Chern connection in the complexified $L_\CC$ of the weight bundle.
Indeed, being $L$ flat, its complexified has a canonical
holomorphic structure given by the $(0,1)$ part of the flat
covariant derivative. Moreover, an adapted  Hermitian metric  can
also be considered, by decreeing the canonical section of $L$ to be
unitary. This way one can look at the Chern connection of $L_\CC$.
A direct computation proves:
\begin{te}\cite{ve}\label{curv}
The curvature of the Chern connection in $L_\CC$ equals
$-2\sqrt{-1}\omega_0$.
\end{te}
This is the result that opens the way for the use of the
algebraic-geometric tools in l.c.K. geometry.

\subsection{Topological properties}
On a Vaisman manifold, one can consider differential operators
adapted to the form $\omega_0$ instead of to the fundamental form
$\omega$. Precisely, let $L_0$ be the wedging operator with
$\omega_0$: $L_0(\eta)=\eta\wedge\omega_0$, let $\Lambda_0$ be its
Hermitian formal adjoint and let $H_0=[L_0,\Lambda_0]$ be their
commutator. If $\xi_1,\ldots,\xi_{n-1}, \sqrt{2}\theta^{1,0}$ is
an orthonormal frame of $T^{1,0}M$, then (cf. \cite{ve}):
\begin{equation}\label{poz}
\begin{split}
\omega_0&=\sqrt{-1}(\xi_1\wedge\bar\xi_1+\cdots+
\xi_{n-1}\wedge\bar\xi_{n-1}),\\
H_0(\eta)&=(p-n+1)\eta,
\end{split}
\end{equation}
for any
$\eta=\xi_{i_1}\wedge\cdots\wedge\xi_{i_k}\wedge\bar\xi_{i_{k+1}}
\wedge\cdots\wedge\bar\xi_p\wedge R$,
where $R$ is a monomial in
$\theta^{1,0}$ and $\theta^{0,1}$ only and $p$ is the number of
$\xi$-s in $\eta$. Verbitsky then shows that, as in the case of
K\"ahler geometry, the triple $L_0, \Lambda_0, H_0$ is a
$\mathrm{SL}(2)$-triple. He then associates to $H_0$ the weight
decomposition $\Lambda^*(M)=\oplus \Lambda_i^*(M)$ such that the
monomial $\eta$ above has weight $p$. This way one has, at the
level $p$, the decomposition: $\Lambda^p(M)=\Lambda^p_{p-2}(M)
\oplus\Lambda^p_{p-1}(M)\oplus\Lambda^p_p(M)$.  Let
$d_0:\Lambda^p_i(M)\longrightarrow \Lambda^{p+1}_{i+1}(M)$ be the
weight 1 component of the usual de Rham differential and denote
with $\partial_0$ and $\bar \partial_0$ its $(1,0)$ and $(0,1)$
part respectively. One immediately checks that $d_0\omega_0=0$,
$d_0\theta=0$, $d_0J\theta=0$. Using algebraic techniques
introduced by Grothendieck (algebraic differential operators),
Verbitsky shows that
$$d_0^2=\partial_0^2=\bar\partial_0^2=0$$
and
he computes the Kodaira identities for the new introduced
operators:
\begin{te}\cite{ve}
On a Vaisman manifolds, the above operators satisfy the following
commutation relations:
\begin{equation}\label{comm}
\begin{split}
[\Lambda_0,\partial_0]&=\sqrt{-1}\bar\partial_0^*, \quad
[L_0,\bar\partial_0]=-\sqrt{-1}\partial_0^*\\
[\Lambda_0,\bar\partial_0^*]&=-\sqrt{-1}\partial_0, \quad
[L_0,\bar\partial_0^*]=\sqrt{-1}\bar\partial_0.
\end{split}
\end{equation}
\end{te}
Here $*$ denotes the Hermitian adjoint.

The same commutation formulas can then be obtained for the similar
operators acting on sections of a holomorphic bundle
$V\longrightarrow M$ endowed with a Chern connection with respect
to which $d_0$ is now the weight 1 component. Preserving the same
notations as above, but understanding now the operators with
coefficients in $V$, one proves again the formulas \eqref{comm}.
This implies the following vanishing result:
\begin{te}\cite{ve}\label{an}
Let $M$ be a compact Vaisman manifold and let $V$ be any positive
tensor power of its weight bundle $L$ equipped with its Hermitian
structure. Let $\bar\partial_0:\Lambda^{0,p}\otimes
V\longrightarrow \Lambda^{0,p+1}\otimes V$ be the
$\omega_0$-differential associated to the Chern connection in $V$
as above. Then all $\bar\partial_0$-harmonic forms are trivial in
dimension $p\leq \dim_\CC M-2$. Moreover, $\eta(\theta^\sharp)=0$ for each
$(0,p)$ $\bar\partial_0$-harmonic form $\eta$ with $p=\dim_\CC-1$.
\end{te}
For the proof, one first shows, as in the case of the operators on
K\"ahler manifolds, using the commutation formulas, the following
connecting formula between the Laplacians of $\partial_0$ and
$\bar\partial_0$:
$$\Delta_{{\partial}_0}-\Delta_{\bar{\partial}_0}=\sqrt{-1}[\Theta_V,
\Lambda_0],$$ where $\Theta_V$ is the curvature operator of the
Chern connection in $V$. As by Theorem \ref{curv}
$\sqrt{-1}\Theta_V=c\omega_0$, $c>0$, one has
$$\Delta_{\bar{\partial}_0}=\Delta_{{\partial}_0}+cH_0.$$
But from the second formula in \eqref{poz}, $H_0$ is positive
definite on $(r,0)$-forms for $r<n-1$. Hence
$\Delta_{\bar{\partial}_0}$ is the sum between a positive definite
and a positive semi-definite operator and thus all his
eigenvalues are strictly positive. This implies that no harmonic
$(r,0)$-forms exist for $r<n-1$. The last statement follows
similarly.

\medskip

On the other hand, if one refers to the notion of basic cohomology
with respect to the canonical foliation $\mathcal{F}$, it is clear
that on basic forms, the operators $\bar\partial$ and
$\bar\partial_0$ coincide. This way, a Tsukada's result about
basic cohomology can be expressed as follows:
\begin{te}\cite{ts}\label{basic}
Let $\eta$ be a $(p,q)$-form, $p+q\leq n-1$, on a compact
$n$-dimensional Vaisman manifold. The following conditions are
equivalent:\\
$(i)$ $\eta$ is $\bar\partial$-harmonic.\\
$(ii)$ $\eta=\theta^{0,1}\wedge\al+\be$, where $\al$, $\be$ are
$\bar\partial_0$-harmonic basic forms and satisfy
$\Lambda_0\al=\Lambda_0\be=0$.
\end{te}

For a first application of these results, suppose that the
canonical bundle $K$ of the Vaisman manifold is a negative power
of $L$. This happens, for example, if $M$ is
Einstein-Weyl\footnote{It is in fact redundant to say Einstein-Weyl
Vaisman manifold, as any Hermitian-Einstein-Weyl manifold is
Vaisman, cf. \cite{pps}.} (cf. \cite{ve}). Now  look at basic
$\bar\partial$-harmonic forms. By the above, these are precisely
the $\bar\partial_0$-harmonic forms. By Serre duality
$\mathcal{H}^p_{\bar\partial_0}(\mathcal{O}_M)\cong
\overline{\mathcal{H}^{n-p}_{\bar\partial_0}(K)}$. But from
Theorem \ref{an},
$\overline{\mathcal{H}^{n-p}_{\bar\partial_0}(K)}=0$ for $p>1$ and
for $p=1$ the non-trivial $\eta\in
\mathcal{H}^{n-p}_{\bar\partial_0}(K)$ satisfy
$\eta({\eta^{0,1}}^\sharp)=0$. Again by Serre duality one derives
that the $\bar \partial_0$-harmonic $(0,p)$-forms in
$\mathcal{H}^p_{\bar\partial_0}(\mathcal{O}_M)$ vanish for $p>1$
and are proportional to $\theta^{0,1}$ if $p=1$. Hence in this
latter case the form cannot be basic, thus there are no
non-trivial basic $\bar\partial$-harmonic $(0,p)$-forms. Summing
up, one obtains the following information about the cohomology of
the structural sheaf:
\begin{te}\cite{ve}\label{ubu}
Let $M$ be a compact Vaisman manifold whose canonical bundle is a
negative power of the weight bundle. Then
$\mathcal{H}^i(\mathcal{O}_M)=0$ for $i>1$ and
$\mathcal{H}^1(\mathcal{O}_M)=1$.
\end{te}

\medskip

In fact, this  result was independently proved by B. Alexandrov and
S. Ivanov  in a more general setting (here I quote  it adapted to
the present situation):

\begin{te}\cite{ai}\label{ubuu}
Let $M$ be a compact Vaisman manifold. Suppose that the Ricci
tensor of its Weyl connection is non-negative. Then $b_1(M)=1$ and
the Hodge numbers satisfy:
\begin{equation*}
\begin{split}
h^{p,q}(M)&=0,\quad p=1,2,\ldots, n\\
h^{0,q}(M)&=0, \quad q=2,3,\ldots, n\\
h^{0,1}(M)&=1.
\end{split}
\end{equation*}
\end{te}
The proof of Alexandrov and Ivanov  is different in spirit from
Verbitsky's: it relies on a Bochner type formula in Weyl geometry.
If $D$ denotes the Weyl covariant derivative on $M$ and $Ric^D$
the symmetric part of its curvature tensor (but note that on
Vaisman manifolds the Ricci tensor of the Weyl connection is
symmetric), they prove the following formula:
\begin{equation}\label{ric}
Ric^D(X,X)\geq \frac{(n-2)(n-4)}{8}\{\mid X\mid^2-\theta(X)^2\}.
\end{equation}
It immediately follows that $b_1(M)\leq 1$ and, as on compact
Vaisman manifolds the first Betti number is odd, it must be $1$. On the
other hand, the Chern covariant derivative $\nabla^C$ on $M$, is
related to the Weyl covariant derivative by $$\nabla^C=D+\frac 12
\theta\otimes Id+ \frac 12 J\theta\otimes J.$$ This gives the
relation between the respective curvature tensors:
$$R^C=R^D+\frac 12 dJ\theta.$$ 
Let
$$R^C(\omega)=\frac 12\sum g(R^C(e_i, Je_i)X,Y),$$
$\{e_i\}$  being an orthonormal frame. Let
$k^C(X,Y)=R^C(\omega)(JX,Y)$ be the associated symmetric,
$J$-invariant bilinear form. One has
$$k^C=Ric^D+\frac 12 \langle
dJ\theta, \omega\rangle g.$$ But it is not hard to prove that
$$dJ\theta=\omega+\theta\wedge J\theta,$$
so that
$$k^C=Ric^D+(n-1)g.$$
As $Ric^D\geq 0$,  $k^C>0$ follows. Now one  applies the classical
vanishing theorems for holomorphic forms and derive $h^{p,0}=0$
for $p=1,2,\ldots,n$. For the rest, one uses the following
relations proved in \cite{ts}:
$$h^{m,0}=h^{m-1,0}, \quad h^{0,p}=h^{p,0}+h^{p-1,0}, \, p\leq
n-1, \quad 2h^{1,0}-b_1-1, \quad 2h^{0,1}=b_1+1$$ and Serre
duality.

\medskip

Another vanishing result I want to mention is the following:
\begin{te}\cite{ve}\label{hol}
Let $M$ be a compact Hermitian-Einstein-Weyl manifold. Then all
holomorphic $p$-forms on $M$ vanish, for $p>0$ and $h^1(M)=1$.
\end{te}
Note that in general $h^1$ can be arbitrarily big (still odd) on a
Vaisman manifold: it is enough to think at the total space of the
induced Hopf bundle over a projective curve of big genus.

The second statement of the theorem follows from the first one by
looking at the Dolbeault spectral sequence $E^{p,q}_r$. Since
$E_2^{p,q}=H^q(\Omega^p(M))$, one gets  $H^0(\Omega^1(M))=0$ and,
from Theorem \ref{ubu}, $H^1(\Omega_0(M))=\CC$. Thus
$$h^1(M)\leq \dim H^0(\Omega^1(M))+\dim H^1(\Omega_0(M)) =1.$$
Now, if $h^1=0$, the weight bundle would have trivial monodromy,
contradiction with it being isomorphic to $\ZZ$.

As for the first statement of the theorem, it follows by first
showing that $TM$, as a holomorphic Hermitian bundle, satisfies
$\Lambda_0(\Theta_{TM})=-c\sqrt{-1}Id_{TM}$, $c>0$ (in other
words, it is $\omega_0$-Yang-Mills with negative constant), then
showing that such a bundle cannot have global holomorphic
sections.

\medskip

 On a l.c.K. manifold one may also consider the twisted
differential operator $d^\theta$ acting as
$$d^\theta \eta = d\eta - \theta\wedge\eta.$$
It is immediate that $\dt\circ\dt=0$, hence $\dt$ produces a
cohomology which groups are denoted with $H^*_\theta(M,\RR)$. It
was first introduced and studied in \cite{gl} for locally
conformally symplectic manifolds. For compact manifolds, the top
$\dt$ cohomology was known (\emph{loc. cit.}) to vanish. However,
it was recently proved that on compact Vaisman manifolds (in fact, the
original proof refers, more generally, to Riemannian manifolds
endowed with a parallel one-form) all these cohomology groups
vanish (see also Remark \ref{loc_syst} for a proof which works only on Vaisman manifolds):
\begin{pr}\cite{llmp}\label{dtcoh}
On a compact, non-K\"ahler Vaisman manifold all the cohomology
groups $H^*_\theta(M,\RR)$ are trivial.
\end{pr}
The proof is technical and is based on the following relations
(here one takes into account that $B=\theta^\sharp$ is Killing etc.):
\begin{equation*}
\Ll_B=-d^*\circ e(\theta)-e(\theta)\circ d^*, \quad d^*\circ
\Ll_B=\Ll_B\circ d^*,
\end{equation*}
where $e(\theta)$ denotes the exterior product operator with
$\theta$ and $d^*$ the Riemannian adjoint of $d$. If  $i(\theta)$
is the interior product with $\theta$ and ${\dt}^*$  the formal
adjoint of $\dt$, then one may  derive the formulas:
\begin{equation*}
\begin{split}
\dt\circ i(\theta)&=-i(\theta)\circ\dt+\Ll_B+Id, \quad
{\dt}^*\circ i(\theta)=-i(\theta)\circ {\dt}^*,\\
\dt\circ \Ll_B&=\Ll_B\circ\dt,\quad \quad \quad \quad \quad \quad
\quad {\dt}^*\circ \Ll_B=\Ll_B\circ{\dt}^*.
\end{split}
\end{equation*}
With this, one may show that $\langle \Ll_B\eta, \eta\rangle=0$
for any $k$-form $\eta$. Hence, if $\dt \eta=0$, then
$\Ll_B\eta=-\eta$ and $\eta=0$ follows.

\begin{re}
Recently, A. Banyaga computed the $d^\theta$ cohomology of a
4-dimensional compact locally conformally symplectic manifold
which is diffeomorphic with the Inoue surface of the first kind
and showed that this $d^\theta$-cohomology is non-trivial, cf.
\cite{ban}. In view of the previous result, the reason might be
that, as Belgun proved (see above), as a complex surface, this
manifold does not admit any metric with parallel Lee form.
\end{re}

\subsection{The structure of compact Vaisman manifolds}
From the proof of Theorem \ref{flux} (see also \cite{gop}) it
follows that the Riemannian cone of a Sasakian manifold admits a
globally conformally K\"ahler metric with parallel Lee form.
Conversely, the universal cover of any Vaisman manifold admits
such a structure, the Sasakian manifold being the universal cover
of leaves of the foliation $\Ker \theta$.

So, to obtain Vaisman manifolds one has the following receipt:
pick a Sasakian manifold $W$, construct the Riemannian cone, then
identify a group of holomorphic homotheties of the K\"ahler metric
on the cone. One natural way of finding such a group is to select
a Sasakian automorphism $\f$ of $W$ and to consider the group
generated on the cone by the transformation $(w,t)\mapsto (\f(w),
qt)$, $q\in \RR$, $q>1$. This is a homothety with scale factor
$q^2$. The quotient $M_{\f,q}$ is the Riemannian suspension of
$\f$ over the circle of length $2\pi q$ and is a Vaisman manifold
(cf. \cite{go}, \cite{ko}, \cite{gop}, \cite{ov1}).

This construction can be reversed, thus providing the general
structure of compact Vaisman manifolds:

\begin{te}\cite{ov1}\label{strte}
A compact Vaisman manifold $M$ admits a canonical Riemannian
submersion $p:M\longrightarrow S^1$ with fibres isometric to a
Sasakian manifold $W$. Moreover, there exists a Sasakian
automorphism $\f$ of $W$ such that $M$ is isomorphic with the
Vaisman manifold $M_{\f,q}$ above.
\end{te}

The proof of the theorem is a straightforward application of
Theorem \ref{mono}. Indeed, let $\{\gamma_1,\ldots,\gamma_k\}$ be
a basis of $H^1(M)$, let $\al_i=\int_{\gamma_i}\theta$ be the
periods of $\theta$ and let $A=\langle \al_1,\ldots,\al_k\rangle$
be the abelian group generated by the periods. One clearly has a
function $p:M\longrightarrow \RR/A$. The monodromy of $L$ being
$\ZZ$, all periods are proportional to a real number $\al$ with
integer coefficient. Hence $A\cong \al\cdot \ZZ$ and one has the
desired map $p$. Passing to the universal cover, one finds a map
$f:\tilde M\longrightarrow S^1$, $df=\theta$. On the other hand,
$\tilde M$ is a cone over a Sasakian $W$. As the monodromy map of
$\tilde M$ preserves both the complex structure and the Lee field
of $\tilde M$ and acts as a homothety, it preserves $W$ and
induces its desired Sasakian automorphism $\f$.

\begin{re}

The same result holds, replacing Sasakian by 3-Sasakian, for
locally conformally hyperk\"ahler manifolds. In fact, compact
l.c.h.K. manifolds are, in particular,  l.c.K. Einstein-Weyl
manifolds. As by Theorems \ref{ubuu} and \ref{ubu} such manifolds
have $b_1=1$, the above construction applies. This was the idea of
the first proof in \cite{ve} for the structure theorem.
\end{re}

One can use Theorem \ref{strte} to determine the fundamental group of a 
compact Vaisman manifold:

\begin{pr}\label{pi1}
The fundamental group of a compact Vaisman manifold is a product 
$K\times \ZZ$ where $K$ is the fundamental group of a compact K\"ahler 
manifold.
 \end{pr}
Indeed, topologically a compact Vaisman manifold is a product of a compact 
Sasakian manifold with a circle. Again topologically, the Sasakian manifold is 
a circle bundle over a compact K\"ahler manifold. Hence the result.

\begin{re}
Note that, in general, nothing can be said about the fundamental group of a l.c.K. manifold, except the fact that it has to send a non-trivial arrow into $\RR_+$.
\end{re}
\begin{re}\label{loc_syst}\cite{ve04a}
Another application of Theorem \ref{strte} is a proof of the vanishing  of 
the $d^\theta$-cohomology on compact Vaisman manifolds, cf. Proposition \ref{dtcoh}. Indeed, this is the cohomology of a locally constant sheaf on $M$ 
corresponding to the weight bundle $L$ considered with its flat connection (a local system). Note that, as topologically $M$ is $W\times S^1$ and the monodromy of $L$ is $\ZZ$, $L$ is the pull-back $p^*L'$ of a local system $L'$ on $S^1$. 
Now, the cohomology of the local system $L$ is the derived direct image $R^iP_*(L)$, 
where $P$ is a projection onto a point. By the above remark and changing the base, $R_iP_*(L)=R^iP_*p^*(L')=R^i({\CC}\otimes L')$, where $\CC$ is viewed as  a trivial local system. From the K\"unneth formula it follows that  $R^i(\CC)$ is a locally constant sheaf on $S^1$, with fiber $H^*(W)$. By the Leray spectral sequence of composition, the hypercohomology of the complex of sheaves $R^*(\CC\otimes L')$ converges to $R^iP_*(L)$. Finally, each $R^i(\CC\otimes L')$ has zero cohomology, hence the spectral sequence vanishes in $E_2$. It then converges to zero.
\end{re}

\begin{re}
Theorem \ref{strte} says that, at least in the compact case, Vaisman
and Sasakian structures are essentially the same. One of the
definitions of a Sasakian structure is the following: a Riemannian
manifold together with a unit Killing $\xi$ field satisfying the
following second order equation: $R(X,Y)\xi=g(\xi, Y)X-g(X,\xi)Y$.
A similar characterization is available also for Vaisman
manifolds:
\begin{pr}\cite{ka}
Let $(M,g)$ be a Riemannian manifold with a unit Killing $1$-form
$\theta'$ and a unit parallel $1$-form $\theta$ orthogonal to
$\theta'$. Then $(M,g)$ admits a compatible Vaisman structure if
and only if:
$$
R(X,Y){\theta'}^\sharp=\theta'(Y)X-\theta'(X)Y-(\theta\wedge\theta')(X,Y)\theta^\sharp.
$$
\end{pr}
\end{re}

The Sasakian manifold exhibited in the structure theorem is ``the smallest'' possible in the following sense:
\begin{pr}\cite{gopp}\label{small}
Let $M$ be a compact Vaisman manifold and $W$ the Sasakian manifold provided by the structure Theorem $\ref{strte}$. Let $W'$ be any Sasakian manifold whose Riemannian cone covers $M$. Then $W'$ covers $W$. In other words, the minimal presentation of a compact Vaisman manifold is $(W, \ZZ)$. In particular, the rank of a compact Vaisman manifold is $1$.
\end{pr}
The argument is as follows: the homotopy sequence associated to the fibration $W\rightarrow M\rightarrow S^1$ gives the exact sequence 
$$0\rightarrow \pi_1(W)\rightarrow \pi_1(M)\rightarrow \ZZ\rightarrow 0.$$
On the other hand, the index of $\pi_1(W')$ cannot be finite in $\pi_1(M)$, as $M$ is compact and the cone $W\times \RR$ non-compact. Hence $\mathrm{coker} \pi_1(W')$, which contains $\ZZ$, is bigger than $\mathrm{coker} \pi_1(W)$. Then, up to conjugation, the image of $\pi_1(W)$ in $\pi_1(M)$ is bigger than the image of $\pi_1(W')$. From here one may conclude that $W'$ covers $W$.

\subsection{Immersing Vaisman manifolds into Hopf manifolds}\label{emb}

The properties of the complexified weight bundle suggest a
Kodaira-Nakano type result for compact Vaisman manifolds. The model
space will now be the Hopf manifold as described in Section
\ref{genex}. The precise statement is:
\begin{te}\cite{ov2}\label{imm}
Let $M$ be a compact Vaisman manifold. Then $M$ admits a holomorphic  immersion
into a diagonal Hopf manifolds $H_\Lambda$, preserving the respective Lee
fields and inducing a finite covering on the image.
\end{te}
Although a better result -- an embedding theorem -- is also available, with a completely different approach, in the framework of l.c.K. manifolds with potential (cf. \S \ref{lck_pot}), I think the proof below is interesting in itself, showing how algebraic geometry methods can be used in this non-algebraic context. I thus present the outline of 
the proof which is divided into two distinct parts: one first assumes $M$
quasi-regular, then shows that any compact Vaisman manifold can be
deformed to a quasi-regular one. Note that the deformation must be
explicitly exhibited, since the l.c.K. class is not closed under small
deformations, cf. \cite{belgun} (see also the discussion in \S \ref{lck_pot} about small deformations of l.c.K. manifolds).

If $M$ is quasi-regular, the push-forward of $L_\CC$ over the base
orbifold $Q$ is a positive line-bundle (the essential technical step
here is to show that $L_\CC$ is trivial along the fibres of $p$). Hence, by the orbifold version
(cf. \cite{Ba}) of the Kodaira-Nakano theorem, $Q$ is projective. One
thus  has an embedding $l:Q\hookrightarrow \CC P^{N-1}$, $N=\dim
H^0(L_\CC^k)$ for an appropriate $k>0$. A basis  of this
space is easily seen to produce functions $\lambda_i:\hhat
M\longrightarrow \CC$ without common zeros on the total space of the weight covering
associated to monodromy group of $L$ (which, remember, is
$\ZZ$). Hence there exists the map $\hhat \lambda:\hhat M \longrightarrow
\CC^N\setminus 0$. This map proves to be compatible with the action of
a generator $\gamma$ of the monodromy representation and with the
multiplication with $q^k$ on $\CC^N$:
\begin{equation*}\label{_tilde_M_immersion_CD_Equation_}
\begin{CD}
\hhat M @>{\hat \lambda}>> \CC^N \setminus 0\\
@V{\gamma}VV  @VV{\text{mult. by $q^k$}}V \\
\hhat M @>{\hat \lambda}>>\CC^N \setminus 0\\
\end{CD}
\end{equation*}
Hence $\hat\lambda$ is the covering of a map $\lambda:M\longrightarrow
(\CC^N \setminus 0)/\langle q^k\rangle$ into a diagonal Hopf
manifold. The reminder of the conclusion follows by observing that
$\lambda$ can be included in the following commutative square:
\begin{equation*}\label{_Lee_fibra_commu_Equation_}
\begin{CD}
M @>{\lambda}>> (\CC^N \setminus 0) \big / \langle q^k\rangle\\
@V{p}VV  @VVV \\
Q @>{l}>> \CC P^{N-1}.
\end{CD}
\end{equation*}

In order to handle the general case, one needs a deformation
argument. The idea is to look at the complexified $G_\CC$ of the group
generated by the Lee flow and at its counter-image $\hhat G_\CC$ by
$\Phi$ acting on the monodromy covering (cf. Subsection
\ref{weight}). The key step consists in proving:
\begin{pr}\label{gc}
$\hhat G_\CC\cong (\CC^*)^l$ for some $l>0$.
\end{pr}
It will be enough to show that $(1)$ $\hhat G_\CC$ is linear and $(2)$ it has a
compact real form. By the Remmert-Morimoto theorem, $\hhat G_\CC$ is a product
$(\CC^*)^l\times \CC^m\times T$ where $T$ is a torus. The torus
component is easily seen to not appear because the K\"ahler metric on
$\hhat M$ has a potential. To prove $(2)$, one
proceeds as follows: Let $\hhat G_K$ be the Lie subgroup of $\hhat
G_\CC$ generated in $G_\CC$ by the group $\hhat G_0$ of
 K\"ahler isometries of $\hhat G$ and the flow
of $J\theta^\sharp$ (which acts by holomorphic K\"ahler
isometries). Clearly $\Phi(\hhat G_K)$ consists of isometries of
$M$. Moreover, one proves that $\Phi:\hhat G_K\rightarrow
\mathrm{Iso}(M)$ is a monomorphism, hence $\hhat G_K$ is compact as a
subgroup of the isometry group.  As $\hhat G\cong (S^1)^k\times \RR$
is generated by $\hhat G_0$ and $e^{t\theta^\sharp}$ and $\hhat G_K$ is
generated by $\hhat G_0$ and $e^{tJ\theta^\sharp}$. Hence
the complexifications of $\hhat G_K$ and $\hhat G$
coincide. All in all one sees that $\hhat G_K$ is a real form of $\hhat
G_\CC$.

\smallskip

The next step is to deform the given Vaisman structure still remaining in the
class of Vaisman structures. To do this, one
picks an element $\gamma'$ sufficiently close in $\hhat G_\CC$ to the
generator $\gamma$ of the monodromy representation
$\ZZ\cong\Gamma\subset \hhat G$. Let $\Gamma'$ be the abelian subgroup
generated in $\hhat G_\CC$ by $\gamma'$. Recall (cf. \cite[Pr. 4.4]{ve}) that the K\"ahler form of $\hhat
M$ is exact and has a potential $\mid V\mid^2$, where $V$ is the lift to $\hhat
M$ of the Lee field. Then the squared length of $V'=\log \gamma'$ will still be
a K\"ahler potential with corresponding K\"ahler form $\omega'$. It can be seen
that the flow of $V'$ acts by holomorphic conformalities with respect to
$\omega'$, hence the quotient $\hhat M/\Gamma'$ is l.c.K. with, by construction,
parallel Lee form. This proves:

\begin{pr}
The quotient $\hhat M/\Gamma'$ is a compact Vaisman manifold with Lee field
proportional to $V'=\log \gamma'\in \mathrm{Lie}(\hhat G_\CC)$\footnote{Note that
by a result in \cite{ts}, the Lee field is unique on a compact manifold of
Vaisman type.}.
\end{pr}

What one  needs is in fact a quasi-regular deformation. This can be achieved by
further restricting the choice of $\gamma'$. As $\hhat G_\CC$ has a compact real
form, its Lie algebra can be given a rational structure by saying that
$\delta\in \mathrm{Lie}(\hhat G_K)$ is rational if $e^{t\delta}$ is compact in
$\hhat G_K$. Now, a complex subgroup $e^{z\delta}$ of $\hhat G_\CC$ is isomorphic
with $\CC^*$ if and only if the line $z\delta$ is rational (it contains a
rational point). One may prove:

\begin{pr}
If the line $\CC\log \gamma'$ is rational, then $\hhat M/\Gamma'$ is a
quasi-regular Vaisman manifold.
\end{pr}

As the elements $\gamma'$ satisfying the above assumption are
dense in $\hhat G_\CC$, one can find a quasi-regular Vaisman
deformation $M'=\hhat M/\langle \gamma'\rangle$  of the initial
Vaisman structure. Applying the first step, we obtain a map
$\hhat\lambda':\hhat M\longrightarrow \CC^N\setminus 0$ with the
same properties as above. But $\gamma$ and $\gamma'$ commute. This
gives rise to the following commutative diagram:

\begin{equation*}
\begin{CD}
\hhat M @>{\gamma}>> \hhat M \\
@V{\hat\lambda'}VV  @VV{\hat\lambda'}V \\
H^0(M', L_\CC^k) @>{\hat\lambda'(\gamma)}>> H^0(M', L_\CC^k).
\end{CD}
\end{equation*}

One now chooses $\gamma'$ sufficiently close to $\gamma$ in order that the
eigenvalues of $\hat\lambda'(\gamma)$ be all $>1$. Thus $H:= \left(H^0(M',
L_\CC^k)\backslash 0\right)/\hat\lambda'(\gamma)$ is a (smooth) complex manifold
and one has the immersion $M=\hhat M/\Gamma \hookrightarrow H$.

\smallskip

The final step is to show that $H$ is indeed a Hopf manifold. This amounts to
say that the operator $\hat\lambda'(\gamma)$ is semi-simple. But this follows
easily from the fact that the vector space $\{f:\hhat M\rightarrow \CC \mid
f\circ\gamma'=q^kf\}$ is finite dimensional. This ends the proof of the immersion
theorem.

For the time being, I can only mention a straightforward application of the immersion theorem to Sasakian geometry, namely:
\begin{te}\cite{ov2}
Any compact Sasakian manifold admits a CR immersion, finite covering on the image, into a weighted sphere.
\end{te} 

\subsection{Stable bundles on Vaisman manifolds}
I quote here a recent result obtained in \cite{ve04} on stability on Vaisman and, in particular, diagonal Hopf manifolds.
\begin{te}\cite{ve04}
Let $M$ be a compact Vaisman manifold of complex dimension at least $3$ 
and $V$ a stable bundle on $M$. Then the curvature $\Theta$ of the Chern connection satisfies $\Theta(v,*)=0$ for every $v$ in the canonical foliation $\mathcal{F}$. In particular, $V$ is equivariant with respect to the complex 
Lie group generated by the Lie flow and this equivariant structure is compatible with the connection.
\end{te}
Here stability is understood with respect to the Gauduchon metric.

The proof is divided into two parts. The first one, computational, 
applies to bundles of degree $0$. It follows, in fact, from a more general statement, concerning a Hermitian bundle endowed with a connection and with a primitive 
closed skew-Hermitian $(1,1)$-form 
$\Theta\in \Lambda^{(1,1)}(M,\RR)\otimes_\RR\mathfrak{u}(V)$ (here primitive means that $\Lambda\Theta=0$, $\Lambda$ being the Hermitian adjoint of the operator of tensoring with the Hermitian form).

To pass to bundles of arbitrary degree, observe first that 
a compact Vaisman manifold admits holomorphic Hermitian bundles of arbitrary degree. Indeed, the complexified weight bundle $L_\CC$ on $M$ has degree 
$\delta:=\int \omega_0\wedge\omega^{n-1}>0$. Now, for any $\la\in \RR$, one associates a holomorphic Hermitian bundle with connection $\nabla_{triv}-i\frac{\la}{\delta}\theta^c$. $L_\la$ has degree $\la$. Moreover, $L_\la$ is equivariant with respect to the Lee flow. 

Now tensoring a stable bundle of arbitrary degree with $L_\la$ for appropriate $\la$ produces a stable bundle of degree $0$.

With a proof which is too involved to be reported here, this result leads to the following
\begin{te}\cite{ve04}
Any stable holomorphic bundle on a diagonal Hopf manifold is filtrable.
\end{te}
  
\subsection{Riemannian properties}
The Lee field of the Tricerri metric is harmonic and minimal. For Vaisman manifolds, more can be said: the anti-Lee field also defines a harmonic map in the unit tangent bundle endowed with the Sasaki metric $g^S$. The precise result is:
\begin{pr}\cite{orv} The following properties hold on a Vaisman manifold $(M,g)$:

$i)$ $J\theta^\sharp$  is a harmonic and
minimal vector field. Moreover, it is a harmonic map from $(M,g)$ into the
unit tangent bundle $(T_1M, g^S)$.

$ii)$ The bivector $\theta^\sharp\wedge J\theta^\sharp$, viewed as a map $ (M,g)\rightarrow (G^{or}_2(M), g^S)$ is harmonic and its image is a minimal submanifold.
\end{pr}
I also mention that, for a compact $M$, computing the Hessian of the volume and energy for $J\theta^\sharp$, one may also prove (cf. \cite{orv}) that this vector field is not stable neither as a harmonic map, nor as a minimal submanifold.

On the other hand, another notion of harmonicity for a foliation $\mathcal{F}$ on a Riemannian manifold $(N,g)$ was introduced an studied in \cite{kt} by the equivalence of the following conditions: 

(1) The canonical projection $\pi$ from $TN$ on the normal bundle of the foliation is harmonic;

(2) The leaves of $\mathcal{F}$ are minimal submanifolds in $(N,g)$.

If $N$ is  compact, oriented, $\mathcal {F}$ is Riemannian and $g$ is bundle-like with respect to it, the above conditions are also equivalent with $\Delta \pi=0$.

Using the second characterization and a Bochner-Weitzenb\"ock type formula, the following result was obtained:

\begin{te}\cite{noda}
 Let $\mathcal{F}$ be a harmonic Riemannian foliation with complex leaves on a compact, oriented l.c.K. manifold. Then $\mathcal{F}$ is stable.
 
 In particular, the canonical foliation of a Vaisman manifold is stable.
 \end{te}
 
Note that this notion of harmonicity implies the one  in \cite{ggv}. Indeed, Ichikawa and Noda observe that if $\mathcal{F}$ is harmonic in the sense of Kamber-Tondeur, then its Gauss map is harmonic (cf. Ruh-Vilms' theorem) and hence the foliation is harmonic also in Vanhecke's sense. On the other hand, as observed above, the canonical foliation of a Vaisman manifold is not stable as a harmonic map, hence the two notions of stability are not the same.   

\medskip

Another way of looking at the Riemannian properties of a Vaisman manifold is to consider its transversal K\"ahlerian structure guaranteed by Proposition \ref{kill}. This is the viewpoint adopted in \cite{ba_dr} where transversally holomorphic maps are studied and the following result is obtained as a corollary of more general properties:

\begin{pr}\cite{ba_dr}
Let $\f:M\rightarrow M'$ be a holomorphic or anti-ho\-lo\-mor\-phic map between compact Vaisman manifolds. Then $\f$ is transversally holomorphic and is an absolute minimum for the energy functional in its foliated homotopy class.
\end{pr}

Of course, here also the energy is considered with respect to  the foliated exterior derivative operator etc.   

In the same paper \cite{ba_dr}, the beginning of a theory of transversal l.c.K. structures is sketched. These appear naturally on manifolds which admit submersions over l.c.K. manifolds. I think that the subject still waits to be developed.
\subsection{L.c.K. structures with parallel anti-Lee form}
Another  natural class of l.c.K. manifolds, curiously neglected until recently, is the one defined by the parallelism of the anti-Lee vector field. Although its geometry is completely different from the Vaisman one, I include it in this section. 

This class was studied in \cite{kashiwada02}. I summarize some of the results therein. 

The condition $\nabla (J\theta)=0$ immediately implies $\delta \theta =-(n-1)$, hence, by Green's theorem, the underlying manifold cannot be compact. But the local geometry of such l.c.K. structures is rather rich: 

\begin{pr}\cite{kashiwada02}
On a l.c.K. manifold with parallel anti-Lee form, the distributions $\mathcal{D}_1$  and $\mathcal{D}_2$ locally generated by $J\theta=0$, respectively $\theta=0$ and $J\theta=0$ are integrable. Their leaves are, respectively: non-compact, totally geodesic real hypersurfaces  with an induced Kenmotsu structure\footnote{A Kenmotsu structure is a metric almost contact structure satisfying the integrability condition $(\nabla_X\varphi)Y=-\eta(Y)\varphi X-g(X,\varphi Y)\xi$. A typical example is the product of a K\"ahler manifold $L$ with the line $\RR$, with $\xi=d/dt$, $\eta=dt$, $g=dt^2+e^{2t}g_L$, $\f=\begin{pmatrix} 0&0\\0&J_L\end{pmatrix}$,  cf. \cite{kenmotsu}.}  and totally umbilical complex hypersurfaces on which the induced structure is K\"ahler. 
\end{pr}
Such l.c.K. structures with parallel $J\theta$ appear naturally on the product of a Kenmotsu manifold with the line. But the  converse of this result is still not proved, nor is the result extended to non-trivial bundles.

\section{Locally conformally K\"ahler manifolds with potential}\label{lck_pot}

A wider, proper subclass  of l.c.K. manifolds, containing strictly  the 
Vaisman one was recently defined in \cite{ov3}. The main reason for introducing it was the search for a class stable to small deformations. Roughly speaking, we deal with l.c.K. manifolds that admit a K\"ahler cover having a \emph{global} K\"ahler potential (and hence non-compact) subject to some further restrictions.
\subsection{Definition and examples}
\begin{de}\cite{ov3}
An \emph{l.c.K.  manifold with  potential} is a manifold
which admits a K\"ahler covering $(\tilde M, \Omega)$, and a positive
smooth function $\f:\; \tilde M \rightarrow \RR_{+}$ 
(the \emph{ l.c.K. potential}) satisfying the following conditions:
\begin{enumerate}
\item  It is proper, \emph{i.e.} its level sets are compact.
\item The monodromy map $\tau$ acts on it by multiplication with a constant   
$\tau (\f)=const \cdot \f$.
\item It is a  K\"ahler potential, \emph{i.e.} 
$\sqrt{-1} \partial\bar\partial \f = \Omega$.
\end{enumerate}
\end{de}

Here the monodromy group $\Gamma$ is simply the deck transformation group of $\tilde M$. If we let 
$\chi:\; \Gamma \rightarrow \RR_{+}$ be the homomorphism mapping 
$\gamma\in \Gamma$ to the number $\frac{\gamma(\varphi)}{\varphi}$, then it can be shown that, eventually choosing another cover with same properties, $\chi$ may be supposed injective.

Note that the property in the definition is inherited by any closed complex submanifold.

A criterion for the first  condition of the definition is contained in the following

\begin{pr}\cite{ov3}\label{crit_prop}
If the above defined character $\chi$ is injective, then the potential 
$\varphi$ is proper
if and only if the image of $\chi$ is discrete in
$\RR_{+}$.
\end{pr}

\begin{ex}
By \cite[Pr. 4.4]{ve} (cf. also  \S \ref{emb} above), all compact Vaisman manifolds do
have a potential. Indeed, on $\tilde M$, we have $\theta=d\f$, where 
$\f= \Omega(\theta^\sharp, J(\theta^\sharp))$. Using 
the parallelism, one shows that  $\Omega=e^{-\f}\tilde\omega$ is a
K\"ahler form with potential $\f$. 
\end{ex}

\begin{ex}
The two  non-compact l.c.K. manifolds constructed recently in \cite{re} are also examples for the above definition. Here is the outline of the construction which is inspired by the one of the Inoue surfaces.

Let $A=(a_{ij})\in \mathrm{SL}(n,\ZZ)$ a matrix with $n$ eigenvalues $\la_i$ subject to the condition $\la_1>\la_2\cdots>\la_{n-1}>1>\la_n$ and let $D=\log A$ (real logarithm). Let $v_j=(v_{ij})^t$, $j=1,\ldots,n$, real eigenvectors associated to $\la_j$, let $P=(v_{ij})$  and let $Q=P^{-1}$. Define a domain $V$ in $\CC^n$ by
$$V=\left\{\sum\al_iv_i\,\mid\, \al_i\in\CC, \mathfrak{Im} \al_i>0, i=1,\ldots,n-1\right\}\cong H^{n-1}\times\CC.$$
Now, for any field $K=\ZZ, \RR,\CC$, Renaud defines the group $G_K$ equal as a set with $K\times K^n$ and endowed with the following multiplication:
$$(z,b)\cdot(z',b')=(z+z',e^{z'D}b+b').$$
Equivalently, in the eigenbasis $\{v_j\}$, this can be written as
$$(z,b_1,\ldots,b_n)\cdot (z', b_1',\ldots,b_n')=(z+z',e^{z'\rho_1}b_1+b_1',\ldots,e^{z'\rho_n}b_n+b_n'),$$
where $\rho_i=\ln \la_i$.  
One observes that the domain $\tilde{\mathcal V}=\CC\times V$ 
in $G_\CC$ is invariant at the right action of $G_\RR$, hence one can define 
$${\mathcal V}:=\tilde{\mathcal V}/G_\ZZ\hookrightarrow G_\CC/G_\ZZ.$$
The above manifold can be equivalently defined as follows: let $\mathrm{exp}$ be the map 
$$\CC^n\ni (z_1,\ldots,z_n)\mapsto (e^{-2i\pi z_1},\ldots, e^{-2i\pi z_n})\in {(\CC^*)}^n$$
 and call ${\mathcal R}:=\mathrm{exp}(V)$; then $\mathcal V$ coincides with the quotient of $\CC\times \mathcal{R}$ by the group generated by the automorphism
$$(z,w_1,\ldots,w_n)\stackrel{\tilde{g}_A}{\mapsto} (z+1, w_1^{a_{11}}w_2^{a_{12}}\cdots w_n^{a_{1n}}, \ldots, w_1^{a_{n1}}w_2^{a_{n2}}\cdots w_n^{a_{nn}}).$$
Letting $g_A$ be the restriction $\tilde g_A$ to $\mathcal{R}$, one obtains also the manifold $\mathcal{R}/\langle g_A\rangle$.

To prove that $\mathcal{V}$ and $\mathcal{R}/\langle g_A\rangle$ are l.c.K., Renaud explicitly constructs l.c.K. potentials on their universal covers. Let first $\f_i:\mathcal{R}\rightarrow \RR$, $\f_i(w)=\sum_jq_{ij}\ln \abs{w_j}.$ These functions are p.s.h., strictly positive for $i\neq n$, have non-vanishing differential and behave well under the action of $g_A$: $\f_i\circ g_A=\la_i\f_i$. With them, one constructs the potential 
$$\f = \sum_{i=1}^{n-1}\f_i^{r_i}+\f_n^2,$$
where the powers $r_i$ are chosen so that
$$\la_1^{r_1}=\la_2^{r_2}=\cdots=\la_{n-1}^{r_{n-1}}=\la_n^2:=\la.$$
It follows that $r_i$ are strictly negative. Obviously $\f$ is p.s.h. and $\f\circ g_A=\la \f$. As the monodromy is, by construction, $\ZZ$, by Proposition \ref{crit_prop} the potential is also proper, so  $\mathcal{R}/\langle g_A\rangle$ is l.c.K. with potential.

A similar construction is done for $\mathcal{V}$. Renaud then shows that her manifolds are not of K\"ahler type. Her proof uses an idea in \cite{loeb} and 
consists in (1) showing that $\tilde{\mathcal V}$ does not admit any strictly p.s.h. function invariant to the action of $G_\RR$, then (2) arriving at a contradiction by restricting to the 2-dimensional case. It seems to me that this proof is very similar in spirit with Belgun's for the non-existence of l.c.K. metrics on the 3-rd type Inoue surface.

It is also remarkable that these are new examples of non-K\"ahler non-compact manifolds.

In complex dimension 2, these examples recover an Inoue-Hirzebruch surface minus its rational curves. 
\end{ex}

\begin{re} 
I believe that Renaud's examples are not of Vaisman type, but I couldn't verify. Note that without compactness, no general criterion can be applied do decide.
\end{re}

\begin{cex}
Although they are defined using a potential on the universal cover, the manifolds constructed in \cite{ot} do not have a discrete deck group, hence, by Proposition \ref{crit_prop}, they are not l.c.K. manifolds with potential.
\end{cex}

\subsection{Stability at small deformations}
Let now $(M,J,g)$ be an l.c.K. manifold with potential and $(M,J')$ be a small deformation of the underlying complex manifold 
$(M,J)$.
Then $\varphi$ is a proper function on $(\tilde M, J')$ satisfying
$\gamma(\varphi) = \chi(\gamma)\varphi$. It is strictly plurisubharmonic 
because a small deformation of strictly plurisubharmonic
function is again strictly plurisubharmonic. Therefore,
$(\tilde M, J')$ is K\"ahler, and 
$\varphi$ is an l.c.K.-potential on $(\tilde M, J')$. This proves 
\begin{te}\cite{ov3}
The class of compact l.c.K.  manifolds with potential is stable under
small deformations.
\end{te}

In particular, any small deformation of a compact Vaisman manifold is still l.c.K. (with potential), but not necessarily Vaisman. This explains \emph{a posteriori} why the construction in
\cite{go} worked for deforming the Vaisman structure of a Hopf
surface of K\"ahler rank 1 to  a non-Vaisman 
l.c.K. structure on the Hopf surface of K\"ahler rank 0 which, moreover, by
\cite{belgun}, does not admit any Vaisman metric.

The above theorem also proves that the new defined class is strict: not all l.c.K. manifolds admit an l.c.K. potential. For
example, one may consider the blow up in a point of a compact
Vaisman manifold. Moreover, the l.c.K. structure of the Inoue surfaces
do not admit potential, since they can be deformed to the non-l.c.K. 
type Inoue surface $S^+_{n;p,q,r,u}$ with $u\in \CC\setminus \RR$
(cf. \cite{belgun}). 

\subsection{Filling the K\"ahler cover and the embedding theorem}
The key step of the following construction is the observation that, leaving apart the surfaces, the  K\"ahler cover $\tilde M$ supporting the global potential can be compactified with one point to a Stein variety. Precisely:

\begin{te}\cite{ov3}
Let  $M$ be a l.c.K manifold with potential, $\dim M \geq 3$,
and let $\tilde M$ be the corresponding covering. Then $\tilde M$
is an open subset of a Stein variety $\tilde M_c$, with at most
one singular point. Moreover, $\tilde M_c \setminus \tilde M$
is just one point.
\end{te}

The restriction on the dimension is essential to can apply a theorem by Rossi-Andreotti-Siu (cf. \cite[Th. 3, p. 245]{rossi} and \cite[Pr. 3.2]{andreotti_siu}) assuring that the set $\tilde M(a)=\{x\in \tilde M\,\mid\,\f(x)>a\}$, which is holomorphically concave, can be filled, thus being an open set in a Stein variety $\tilde M_c$ with at most isolated singularities. This embedding is then extended to $\tilde M\hookrightarrow \tilde M_c$. To show that  $\tilde M_c$ is indeed obtained from $\tilde M$ by adding just one point, one looks at the generator $\gamma$ of the monodromy $\Gamma\cong\ZZ$ which acts by conformal transformations on $\tilde M_c$. Now, $\gamma$ can be a contraction or an expansion, let us make the first choice: $\gamma(\Omega)=\la\cdot\Omega$, $\la<1$. 
Letting $f$ be any holomorphic function on $\tilde M_c$, it is possible to show that the sequence $\gamma^nf$ converges to a constant. Now, the set $Z:= \tilde M_c \setminus \tilde M$ is compact and is fixed by $\Gamma$. Hence $\sup_Z |\gamma^n f| = \sup_Z |f|$
and $\inf_Z |\gamma^n f| = \inf_Z |f|$. Thus  $\inf_Z |f|=\sup_Z |f|$ for any holomorphic function $f$. Therefore, all  holomorphic functions have the same value in all the points of $Z$. But $\tilde M_c$ is Stein, and thus for any distinct $x$ and $y$ there exists a holomorphic $f$ such that $f(x)\neq f(y)$. This proves that $Z=\{z\}$, one point.

The next step is to show that, in the above hypothesis,

\begin{pr}\cite{ov3}
$\gamma$ acts with eigenvalues strictly smaller than $1$ on the cotangent space $T_z^*\tilde M_c$.
\end{pr}

The proof goes as follows (for the tangent space, which is enough). 
For the smooth case, the argument is the Schwarz lemma. In general, 
one uses the fact that $\gamma^nf=const.$ for any holomorphic $f$ (see above). 
Then, if $\gamma(x)=x$ and $d_x\gamma(v)=\la v$ and $d_xf(v)\neq 0$, 
then $d_x(\gamma^nf)(v)=\la^nd_xf(v)\rightarrow 0$ because $\gamma^nf$ converges to a constant. Hence $\abs{\la}<1$.

This implies that the formal logarithm of $\gamma$ converges. A theorem in \cite{ve_unp} assures that in these conditions $\gamma$ acts with finite Jordan blocks on the formal completion $\widehat{{\mathcal O}}_z$ of $\mathcal{O}_z$ (the local ring of analytic functions in $z\in \tilde M_c$). But one can show that on a Stein variety $S$, for a holomorphic flow with eigenvalues smaller than $1$ on $T_sS$ for some $s$, there exists a sequence of $V_n\subset \mathcal{O}_sS$ of finite dimensional subspaces such that the $s$-adic completion of $\oplus V_n$ is exactly the completion of  $\mathcal{O}_sS$, each $V_n$ being preserved by the flow which acts by linear transformations on it. For the proof, one first observes that $S$ can be assumed smooth, in fact an open ball in  $\CC^n$: $S$ can be seen as an analytic subvariety of an open ball and  the holomorphic flow can be extended to that ball. Now it is possible to apply an old theorem of Poincar\'e and Dulac (cf. \cite[p. 181]{ar}) which gives the normal form of such a flow: $\la_ix_i+P(x_{i+1},\ldots,x_n)$, where $P$ are resonant polynomials corresponding to the eigenvalues $\la_i$. Each $V_{\la_i}$  generated by the coordinate function $x_i$ is preserved by the flow. Hence, picking some $q_1\in V_{\la_{i_1}}, \ldots, q_k\in V_{\la_{i_k}}$, the flow will preserve the space $V_{\la_{i_1}\ldots\la_{i_k}}$ they generate. On the other hand, the completion of all these $V_{\la_{i_1}\ldots\la_{i_k}}$ is precisely $\widehat{\mathcal O}_sS$.

All in all, we find that it is possible to choose enough holomorphic functions in the ei\-gen\-spa\-ces $V_n$ which, together, give an embedding of the cover $\tilde M_c$ in a $\CC^N$ in such a way that the monodromy $\Gamma$ acts equivariantly on $\CC^N$, with eigenvalues smaller than $1$. Hence $M$ embeds in $(\CC^N\setminus 0)/\Gamma$. Such a quotient was called a linear \emph{Hopf manifold}, generalizing both class  0 and class 1 Hopf surfaces of Kodaira\footnote{The terminology is misleading, but belongs to Kodaira. In fact, the group that defines the class 0 Hopf surfaces does not contain  linear transformation except when the transformation is a Jordan cell.}. As in \cite{go}, such a Hopf manifold is l.c.K., but non-Vaisman in general.

\begin{re}
1) The Hopf manifolds just constructed are the appropriate generalization to arbitrary dimension of the Hopf surfaces of class 0.

2) Using another approach, namely exploiting Sternberg's normal form of holomorphic contractions, Belgun also showed that Hopf manifolds are l.c.K. and classified the Vaisman ones among them: they correspond to contractions whose normal form is resonance free, see (and mainly listen to) \cite{be02}.
\end{re}

Summing up, we arrive at the following result which improves a lot the immersion Theorem \ref{imm}:

\begin{te}\cite{ov3} Any compact l.c.K. manifold with potential, of complex dimension at least $3$, admits an embedding in a linear Hopf manifold. If $M$ is Vaisman, it can be embedded in a Vaisman-type Hopf manifold $(\mathbb{C}^N\setminus 0)/\langle A\rangle$, where A is a diagonal linear operator with eigenvalues strictly smaller than $1$.
\end{te} 

The first statement was already proved. When $M$ is Vaisman, the idea is to first identify the $L^2$-completion of the space $V$ of CR-holomorphic fuctions on $S$ with the space of 
$L^2$-integrable holomorphic functions on a pseudoconvex domain bounded by a level set $S_a$ of the potential function. This is assured by a result in \cite{md}. One next considers the flow $X$ of the Lee field, which is holomorphic 
and naturally acts on $V$. As $X$ acts by homotheties on any K\"ahler cover, it acts on $V$ as a self-dual operator, hence it is diagonal any finite-dimensional eigenspace. In particular, this holds for the subspaces 
$V_{\la_{i_1}\ldots\la_{i_k}}$ constructed above which were used to construct the embedding of $\tilde M_c$ in $\CC^N$.

As a straightforward consequence, one obtains the embedding result for Sasakian manifolds:

\begin{te}\cite{ov3} Any compact Sasakian manifold admits a CR-em\-bed\-ding into a Sasakian weighted sphere, preserving the respective Reeb fields.
\end{te}

\section{Locally conformally K\"ahler reduction}\label{red}
Symplectic reduction (at $0$) was easily extended to the K\"ahlerian case.
In fact, it was enough to verify that if the Hamiltonian action
of a group preserves  also the complex structure, then this is
projectable on the symplectic quotient and still compatible with
the symplectic form. But the passage from holomorphic isometric
actions to holomorphic \emph{conformal} actions was more
difficult. Among the main problems one has to solve are: 

(1)
producing a good definition of specific action on l.c.K.
manifolds, producing a momentum map; 

(2) producing a reduction
scheme compatible with the reduction of the related structures
(K\"ahler, Sasakian); 

(3) finding conditions for the quotient of a
Vaisman manifold to be Vaisman. 

\noindent One sees that a major difficulty
is that on l.c.K. manifolds one should act with conformalities,
whereas for the other structures one uses isometries.

The first results were obtained in \cite{bg}. This paper
essentially contains all that is needed to define and perform
l.c.K. reduction, the theory being sketched in the language of
conformal geometry. The symplectic version of this result, namely
the locally conformally symplectic reduction, was independently
discovered  in \cite{hr} and presented in a
local language. Then, in \cite{gop}, locally conformally
symplectic reduction is shown to be compatible with the complex
frame and, on the other hand, the reduction thus obtained is shown
to be equivalent with the Biquard-Gauduchon one. Moreover, all the
above problems are addressed and  solved. Further developments are given in \cite{gopp}.

To begin with, using the twisted differential $d^\theta$, one
defines a twisted Poisson bracket by
$\{f_1,f_2\}^\theta=\omega(\sharp d^\theta f_1, \sharp d^\theta
f_2)$ which is easily seen to satisfy the Jacobi identity. Now an 
action of a subgroup $G\subseteq \mathrm{Aut}(M)$ is said
\emph{twisted Hamiltonian} if there exists a Lie algebra
homomorphism (with respect to the twisted Poisson bracket)
$\mu^\cdot:\mathfrak{g}\rightarrow \mathcal{C}^\infty(M)$ such
that $i(X_M)\omega=d^\theta\mu^{X_M}$ for any vector field $X\in
\mathfrak{g}$\footnote{I denote with $X_M$ the fundamental field associated to 
$X$, namely whose flow is $\frac{d}{dt}e^{tX}\cdot x$.}. 

An essential property of a twisted Hamiltonian action is its conformal invariance in the following sense: by direct computation one sees that to the conformal change $\omega'=e^\al\omega$ correspond the following ``conformal'' changes:
\begin{equation*}
\sharp_{\omega'}d^{\theta'}(e^\al f)=\sharp_\omega d^\theta f,\quad 
\{e^\al f_1,e^\al f_2\}^{\omega'}=e^\al\{f_1,f_2\}^\omega.
\end{equation*}
This implies that multiplication with $e^\al$ induces an isomorphism of Lie algebras between $(\mathcal{C}^\infty(M), \{,\}^\omega)$ and  $(\mathcal{C}^\infty(M), \{,\}^{\omega'})$ which commutes with taking Hamiltonians.

When an action is twisted Hamiltonian,  $\mu$ is called a momentum map for the
action of $G$. It is important to note that the definition is
consistent with the conformal framework: if $g'=e^\alpha g$, then
$\mu'=e^\alpha\mu$ and, in particular, the level set at $0$ is
well-defined. That is why, as in the classical contact (and
Sasakian) case, one performs, for the moment at least, only l.c.K.
reduction \emph{at zero}.

 The l.c.K. reduction can now be stated in:
\begin{te}\cite{gop}
Let $(M,J,g)$ be a locally conformally K\"ahler manifold and $G$ a
subgroup of $\mathrm{Aut}(M)$ whose action is twisted Hamiltonian.
Suppose that $0$ is a regular value for the associated momentum
map $\mu$ and that the action of $G$ is free and proper on
$\mu^{-1}(0)$. Then there exists a locally conformally K\"ahler
structure $(J_0,g_0)$ on $M_0=\mu^{-1}(0)/G$, uniquely determined
by the condition $\pi^*g_0=i^*g$ where
$i:\mu^{-1}(0)\rightarrow M_0$ is the canonical projection.
\end{te}

The striking property, cf. \cite{gop}, of twisted Hamiltonian
actions is that they lift to Hamiltonian (in particular isometric)
actions with respect to the K\"ahler metric of any globally
conformally K\"ahler covering: indeed, the lifted action is
twisted Hamiltonian with respect to the lifted l.c.K. metric; but this one is globally conformal to a K\"ahler metric with respect to which the lifted action is still twisted Hamiltonian (see the above remark on  the conformal invariance). As the Lee form of a K\"ahler metric is zero,  the
notion of twisted Hamiltonian coincides with that of Hamiltonian
action. This reduces l.c.K. reduction to K\"ahler reduction of any
K\"ahler covering:

\begin{te}\cite{gop}
Let $(M,J,g)$ be a locally conformally K\"ahler manifold and let
$G\subset \mathrm{Aut}(M)$ satisfy the hypothesis of the above
reduction theorem. Let  $\tilde G$ be a lift of $G$ to the
universal covering $\tilde M$ of $M$. Then the K\"ahler reduction
is defined, with momentum map  denoted by $\mu_{\tilde M}$,
$\tilde G$ commutes with the action of $\pi_1(M)$, and the
following equality of locally conformally K\"ahler structures holds:
\begin{equation}\label{ugualikaehler}
\mu^{-1}(0)/G\cong(\mu_{\tilde M}^{-1}(0)/\tilde G)/\pi_1(M).
\end{equation}

Conversely, let $\tilde G$ be a subgroup of isometries of a
K\"ahler manifold $(\tilde M, \tilde{g},J)$ of complex dimension
bigger than $1$ satisfying the hypothesis of K\"ahler reduction and
commuting with the action of a subgroup
$\Gamma\subset\mathcal{H}(\tilde M)$ of holomorphic homotheties
acting freely and properly discontinuously and such that
$\rho(\Gamma)\neq 1$. Moreover, assume that $\Gamma$ acts freely
and properly discontinuously on $\mu_{\tilde M}^{-1}(0)$. Then
$\tilde G$ induces a subgroup $G$ of $\mathrm{Aut}(M)$, $M$ being
the locally conformally K\"ahler manifold $\tilde M/\Gamma$, which
satisfies the hypothesis of the reduction theorem, and the
isomorphism~\eqref{ugualikaehler} holds.
\end{te}

This construction applies in particular to the Riemannian  cone
which covers a (compact) Vaisman manifold and permits the link
between l.c.K. and Sasakian reductions (for the Sasakian reduction,
see \cite{gro}). The precise result reads:

\begin{te}\cite{gop}\label{compatibilitasasaki}
Let $W$ be a Sasakian manifold and let $\Gamma$ be  a group of
Sasakian automorphisms inducing holomorphic homotheties on the
cone $W\times \RR$. Let $M=W\times \RR/\Gamma$ be the
corresponding Vaisman manifold.
 Let
$G\subset \mathrm{Iso}(W)$ be a subgroup  satisfying the
hypothesis of Sasakian reduction. Then $G$ can be considered as a
subgroup of $\mathcal{H}(W\times\RR)$. Assume that the action of
$G$ commutes with that of $\Gamma$, and that $\Gamma$ acts freely
and properly discontinuously on the K\"ahler cone
$(\mu_W^{-1}(0)/G)\times\RR$.

Then $G$ induces a subgroup of $\mathrm{Aut}(M)$ satisfying the
hypothesis of the reduction theorem, and the reduced manifold is
isomorphic with $((\mu_W^{-1}(0)/G)\times\RR)/\Gamma$. In
particular the reduced manifold is Vaisman.
\end{te}

In the compact case, this situation can be, rather unexpectedly, reversed, hence l.c.K. reduction of compact Vaisman manifolds is  completely equivalent to the Sasakian reduction. In general, if one starts with a K\"ahler action on the cone over a Sasakian manifold, this action does not come from a Sasakian one, namely the action does not necessarily commute with the translations along the generators of the cone. But if the Sasakian manifold is compact, then:

\begin{pr}\cite{gopp}
Any K\"ahler automorphism of the cone $W\times \RR$ over a compact Sasakian manifold  $W$ is of the form $(f,Id)$, where $f$ is a Sasakian automorphism 
of $W$.
\end{pr}

The proof is based on showing that, as a metric space -- with the distance induced  from the cone Riemannian metric -- the cone can be completed with only one point. This allows extending any isometry to the completion by fixing the new added point which, in turn, implies that the isometry preserves the horizontal sections of the cone, so it comes from an isometry $f$ of $W$. The fact that $f$ preserves the Sasakian structure is an easy consequence of the K\"ahlerian character of the initial automorphism.

\begin{re}
Due to the fact that the above proof is, in its essence, Riemannian, the result also works for almost K\"ahler cones over compact $K$-contact manifolds and, in general, for the Riemannian cone over a compact contact metric manifold.
\end{re}

As the universal cover of a compact Vaisman manifold is a Riemannian  cone over a compact Sasakian manifold, due to the above proven compatibility of reductions, one obtains:
\begin{te}\cite{gopp}
The l.c.K. reduction of a compact Vaisman manifold is a Vaisman manifold.
\end{te}
 
This mechanism allows the obtaining of examples of Vaisman
manifolds by merely reducing weighted Sasakian spheres (cf. \cite{gop}). 
In fact, let $S^1$ act on the weighted sphere $S^{2n-1}$ endowed with the contact form $\eta_A$ and the metric described in paragraph \ref{genex} by 
\begin{equation}\label{s1}
(z_1,\ldots,z_n)\mapsto (e^{i\al_1t}z_1,\ldots,e^{i\al_nt}z_n).
\end{equation}
This action is by Sasakian automorphisms, independently on the weights $A=(a_1,\ldots,a_n)$. The corresponding Sasakian momentum map reads:
$$\mu(z)=\frac{1}{2(\sum \al_i\mid z_i\mid^2)}\sum \al_i\mid z_i\mid^2.$$
Its zero level set is non-empty as soon as the coefficients $\al_i$ have not the same sign. Hence, it is no loss of generality to suppose the first $k$ $\al$-s negative and the others positive. Then $\mu^{-1}(0)$ is diffeomorphic with $S^{2k-1}\times S^{2n-2k-1}$. To have a good quotient, the action of the circle on $\mu^{-1}(0)$ has to be free and proper. A simple computation shows that the necessary and sufficient conditions is that all the positive $\al$-s be coprime with the negative ones. Applying the Sasakian reduction as in \cite{gro}, this provides Sasakian quotients, call them $S(\al)$, of the weighted sphere. So, the above discussion yields:

\begin{pr}\cite{gop}
For any $(\al_1,\ldots,\al_n)$ as above, for any $(a_1,\ldots,a_n)\in\mathbb{R}^n$, $0<a_1\leq\cdots\leq a_n$, for any $(c_1,\ldots,c_n)\in ({S^1})^n$, there exists a Vaisman structure on $S(\al)\times S^1$ which is the l.c.K. reduction of the Vaisman structure of the Hopf manifold associated to the couple $((a_1,\ldots,a_n),(c_1,\ldots,c_n))$ with respect to the $S^1$ action \eqref{s1}.
\end{pr} 

For example, if $n=2k=4s$, then $\mu^{-1}(0)$ is diffeomorphic to $S^{4s-1}\times S^{4s-1}$ and the Sasakian quotient can be seen to be diffeomorphic to $S^{4s-1}\times \CC\mathbb{P}^{2s-1}$. This leads to Vaisman structures on  $S^{4s-1}\times \CC\mathbb{P}^{2s-1}\times S^1$.

Let me stress that when reducing a compact Vaisman manifold $M$, one in fact reduces the simply connected Sasakian manifold whose cone is the universal cover of  $M$. But this Sasakian reduction need not be simply connected, hence the cone over it is a cover, but not necessarily the universal cover of the reduced Vaisman manifold. By contrast, as the Sasakian manifold which is the fibre of the fibration $M\rightarrow S^1$ is the smallest whose cone covers $M$ (see Proposition \ref{small}):
\begin{pr}\cite{gopp}
The reduction procedure is compatible with the structure Theorem $\ref{strte}$.
\end{pr}

But something more general is true: the minimal presentation of a l.c.K. manifold is compatible with the reduction. Indeed one has:

\begin{pr}\cite{gopp}
The minimal presentation of a l.c.K. reduction $M\rid G$ is given by 
\[
(\Kmin\rid G_\text{min}^\circ,\Gammamin).
\]
In particular,the rank is preserved under reduction:
\[
r_{M\rid G}=r_M.
\]
\end{pr}

\section{Open problems}
I shall end these notes with mentioning some open questions that I find worth  thinking about.
\begin{enumerate}
\item The following 3 questions where asked in \cite{ov2} in connection with the algebraic geometry of compact Vaisman manifolds:
\begin{enumerate}
\item Determine the moduli spaces of Vaisman structures. One should note that complex analytic deformations do not preserve the Vaisman class.
\item Given a singular Vaisman variety (in view of the embedding theorem, this is a sub-variety of a diagonal Hopf manifold), does there exist a resolution of singularities within the Vaisman category? Determine the birational maps of Vaisman varieties.
\item Let $M^n$ be a compact Vaisman manifold with canonical bundle isomorphic to $L^{\otimes n}$, where $L$ is the weight bundle (this happens when $M$ is Einstein-Weyl). Does $M$ carry an Einstein-Weyl structure? A positive answer will give a Vaisman analogue of Calabi-Yau theorem.
\end{enumerate}
\item In connection with the last question and with the embedding theorem: Suppose that a compact Hermite-Einstein-Weyl manifold (in particular Vaisman) is isometrically  immersed into a Hopf manifold. What can be said about the Weyl-Ricci curvature? A recent result concerning this problem in the K\"ahler case can be found in \cite{hulin}.
\item  As regards the non-K\"ahler compact surfaces which do not appear in Belgun's list, the most important question is: Which compact complex surfaces with non-zero Euler-Poincar\'e
characteristic do admit l.c.K. metrics?
\item Among the known examples of l.c.K. manifolds, some have discrete fundamental group, some have non-abelian one. The fundamental group of a compact Vaisman manifold is determined in Proposition \ref{pi1}. But in general, which are the groups that can be fundamental groups for a (compact) l.c.K. manifold?
\item Related to the reduction scheme: is it possible to obtain a Vaisman quotient out of a non-Vaisman manifold?
\item As the (compact) Vaisman class is preserved by l.c.K. reduction and Vaisman manifolds are l.c.K. with potential, is the (compact) l.c.K. with potential class preserved by l.c.K. reduction? 
\item Is it possible to construct a convexity theory for l.c.s. and l.c.K. manifolds? The difficulty is that one cannot use Morse theory using the $d$ operator, whereas the cohomology of $d^\theta$ is trivial on compact Vaisman manifolds.
\item What is the structure of compact toric l.c.K. and, in particular, Vaisman manifolds?
\item Not in the very framework of l.c.K. geometry, but still related to it: what is the subclass of locally conformally symplectic manifolds corresponding to Vaisman ones? A natural analogy would be to consider those l.c.s. manifolds which universal cover is globally conformal to a symplectic cone over a contact manifold. Does there exist an intrinsic characterization of this class, eventually in terms of transformation groups as in Theorem \ref{flux}?
\item Find a spinorial characterization of l.c.K. (or, at least, Vaisman) manifolds. On one hand, one can try to characterize l.c.K. manifolds in the Hermitian class as the limiting case of an estimate for the eigenvalues of the Dirac operator. Maybe more at hand would be to start with a l.c.K. manifold endowed with the Gauduchon metric, to find a good notion of a ``twisted Killing spinor'', to obtain a Hijazi type inequality whose limiting case to force the Lee form to be parallel. Such an inequality should refine the one obtained in \cite{mo}, where the complex structure plays no role.
\end{enumerate}

\end{document}